\documentclass[10pt,a4paper]{article}%
\usepackage{ifthen}
\usepackage{amsmath}
\usepackage{amscd}
\usepackage{amssymb}
\usepackage{theorem}
\usepackage[hypertex]{hyperref}%
\usepackage{makeidx}%
\usepackage{latexsym}

\theoremstyle{plain} 
\newtheorem{thm}{Theorem}[section]

\newtheorem{lem}[thm]{Lemma}
\newtheorem{prop}[thm]{Proposition}

\theoremstyle{remark}

\newcommand{\thmref}[1]{Theorem~\ref{#1}}
\newcommand{\propref}[1]{Proposition~\ref{#1}}
\newcommand{\secref}[1]{\S\ref{#1}}

\newcommand{\R}{\mathbb{R}}
\newcommand{\C}{\mathbb{C}}

\newcommand{\pdo}{\psi{\rm do}}
\newcommand{\pdos}{\psi{\rm dos}}

\newcommand{\pdbo}{\psi{\rm dbo}}
\newcommand{\pdbos}{\psi{\rm dbos}}

\newcommand{\Ci}{C^\infty}

\newcommand{\pd}{\partial}

\newcommand{\na}{\nabla}

\newcommand\si{^{-s}}

\newcommand{\noi}{\noindent}

\newcommand{\ol}{\overline}

\newcommand{\too}{\longrightarrow}
\newcommand{\mto}{\mapsto}
\newcommand{\mtoo}{\longmapsto}

\newcommand{\wh}{\widehat}

\newcommand{\wt}{\widetilde}

\newcommand{\bsh}{\backslash}

\newcommand{\sbp}{\mbox{\small{$\partial \pi$}}}

\newcommand{\sbx}{\mbox{\tiny{$\partial X$}}}

\newcommand{\Eb}{E_{\,\mbox{\tiny{$\partial M$}}}}

\newcommand{\tNm}{ \wt{\nabla}}

\newcommand{\Nb}{ \nabla^{\mbox{{\tiny $\partial M$}}}}
\newcommand{\Nm}{ \nabla^{\mbox{{\tiny $M$}}}}
\newcommand{\tNb}{ \wt{\nabla}^{\mbox{{\tiny $\partial M$}}}}

\newcommand{\Nbw}{\nabla^{\,\mbox{{\tiny ${\mathcal W}$}}}}

\newcommand{\Nbk}{\nabla^{\,\mbox{{\tiny ${\mathcal K}$}}}}

\newcommand{\Nsp}{\nabla^{\,\mbox{{\tiny $\sp$}}}}

\newcommand{\Nbkw}{\nabla^{\,\mbox{{\tiny ${\mathcal K},{\mathcal W}$}}}}

\newcommand{\Np}{ \nabla^{\,\mbox{{\tiny $\textsf{P}$}}}}

\newcommand{\tNp}{ \wt{\nabla}^{\,\mbox{{\tiny $\textsf{P}$}}}}
\newcommand{\Nk}{ \nabla^{\,\mbox{{\tiny $P(\textsf{D})$}}}}

\newcommand{\Nzp}{ \nabla^{\,\mbox{{\tiny $\zeta,\textsf{P}$}}}}

\newcommand{\nablab}{ \nabla^{\,\mbox{{\tiny $\partial M$}}}}

\newcommand{\rkw}{R^{\,\mbox{{\tiny ${\mathcal K}, {\mathcal W}$}}}}

\newcommand{\Rb}{ R^{\,\mbox{{\tiny $\partial M$}}}}

\newcommand{\Dsb}{ \Ds_{\,\mbox{{\tiny $\partial M$}}}}

\newcommand{\db}{d{\hskip-1pt\bar{}}\hskip1pt}

\newcommand{\Aa}{{\mathcal A}}

\newcommand{\Cc}{{\mathcal C}}
\newcommand{\Dd}{{\mathcal D}}
\newcommand{\Ee}{{\mathcal E}}

\newcommand{\Hh}{{\mathcal H}}
\newcommand{\Ii}{{\mathcal I}}

\newcommand{\Kk}{{\mathcal K}}

\newcommand{\Mm}{{\mathcal M}}

\newcommand{\Ss}{{\mathcal S}}

\newcommand{\Uu}{{\mathcal U}}

\newcommand{\Ww}{{\mathcal W}}

\newcommand{\as}{{\textsf{a}}}

\newcommand{\gs}{{\textsf{g}}}
\newcommand{\ps}{{\textsf{p}}}
\newcommand{\ws}{{\textsf{w}}}

\newcommand{\As}{{\textsf{A}}}
\newcommand{\Bs}{{\textsf{B}}}
\newcommand{\Cs}{{\textsf{C}}}
\newcommand{\D}{{\textsf{D}}}
\newcommand{\Ds}{{\textsf{D}}}

\newcommand{\Gs}{{\textsf{G}}}

\newcommand{\Is}{{\textsf{I}}}

\newcommand{\Ks}{{\textsf{K}}}

\newcommand{\Ps}{{\textsf{P}}}
\newcommand{\Qs}{{\textsf{Q}}}

\newcommand{\Sss}{{\textsf{S}}}

\newcommand{\Ts}{{\textsf{T}}}

\newcommand{\Ws}{{\textsf{W}}}

\newcommand{\cc}{{\rm c}}

\newcommand{\re}{{\rm Re}}

\newcommand{\<}{\subset}
\newcommand{\ii}{^{-1}}

\newcommand{\G}{\Gamma}

\newcommand{\z}{\zeta}

\newcommand{\s}{\sigma}
\newcommand{\alp}{\alpha}
\newcommand{\bet}{\beta}
\newcommand{\gam}{\gamma}

\newcommand{\lam}{\lambda}

\newcommand{\del}{\delta}

\newcommand{\e}{\varepsilon}
\newcommand{\Gam}{\Gamma}

\newcommand{\Nline}{\mathbb{N}}

\newcommand{\oo}{\infty}

\newcommand{\End}{\mbox{\rm End}\,}

\newcommand{\ch}{\mbox{\rm ch}}

\newcommand{\Gr}{\mbox{${\rm Gr}$}}

\newcommand{\Hom}{\mbox{\rm Hom}}

\newcommand{\ind}{\mbox{\rm ind\,}}

\newcommand{\Ind}{\mbox{\rm Ind\,}}

\newcommand{\Diff}{\mbox{\rm Diff\,}}

\newcommand{\Ker}{\mbox{\rm Ker}\,}
\newcommand{\Cok}{\mbox{\rm Cok}\,}

\newcommand{\tr}{\mbox{\rm tr\,}}
\newcommand{\Tr}{\mbox{\rm Tr\,}}

\newcommand{\g}{\gamma}
\renewcommand{\Re}{{\rm Re}}

\newcommand\Det{\mbox{\rm Det\,}}

\renewcommand{\s}{\sigma}

\newcommand{\Pss}{\Ps_{\mbox{{\tiny $\sigma$}}}}

\newcommand{\ran}{{\rm ran}}

\renewcommand{\sp}{\Sss(\Ps)}

\numberwithin{equation}{subsection}

\begin{document}

\title{$\eta$ Forms and Determinant Lines}

\author{Simon Scott}

\date{}

\maketitle


\section{Introduction}

The purpose here is to give a direct computation of the zeta-function curvature for the determinant
line bundle of a family of APS-type boundary value problems.

Here is the sort of computation we have in mind.

\subsubsection{Example: $\zeta$-curvature on ${\rm CP}^1$}\label{example complex projective space}

Consider the simplest case, $D = id/dx$ over $[\,0,2\pi]$  with
Laplacian $\Delta = -\,d^2/dx^2$. Global boundary conditions for
$D$ are parameterized by $\C P^1$. Specifically, over the dense
open subset of $\C P^1$ parameterizing complex lines $l_z\<\C^2$
given by the homogeneous coordinates $[1,z]$ for $z\in\C$ the
orthogonal projection $P_{z} = \frac{1}{1 + |z|^2}
\begin{pmatrix}
     1 & \ol{z} \\
      z & |z|^2
\end{pmatrix}$
onto $l_z$ parametrizes the boundary condition $P_{z}
\begin{pmatrix}
     \psi(0) \\
      \psi(2\pi)
\end{pmatrix} = 0$; that is,  $\psi(0) = -\ol{z}\,\psi(2\pi)$. Let $D_{P_z}$ denote $D$ with domain restricted to functions satisfying this boundary condition. The adjoint boundary problem is $D_{P^{*}_z}$
with projection  $P^{*}_{z} =  \frac{1}{1 + |z|^2}
\begin{pmatrix}
|z|^2 & -\ol{z} \\
-z &    -1
\end{pmatrix}$
corresponding to  $-z\phi(0) = \phi(2\pi)$. Then $\Delta_{P_z}$ has discrete
  spectrum\\
  $\{(n+\alp)^2, (n-\alp)^2:n\in \mathbb{N}\}$, where $u =
  e^{2\pi i \alp}$ satisfies
 $u^2(1 + |z|^2) + 2u(z + \ol{z}) + (1 + |z|^2) = 0.$  The zeta determinant of $\Delta_{P_z}$ is therefore
\begin{equation}\label{z metric cp1}
{\rm det}_{\zeta}\Delta_{P_z} = 4\sin^2\pi\alp  =
\frac{2|1+\ol{z}|^2}{1 + |z|^2}, \hskip 10mm [1,z]\in\C P^1.
\end{equation}
The Quillen metric evaluated on the holomorphic section identified
with the abstract determinant   $z\mto {\rm det}D_{P_z}$ is
$\|{\rm det}D_{P_z}\|^2 = {\rm det}_{\z}\Delta_{P_z}$ and hence
the canonical curvature $(1,1)$ -form of the determinant line
bundle is
\begin{equation}\label{curv=kahler}
\overline{\pd}\pd\log {\rm det}_{\z}\Delta_{P_z} = \frac{dz\wedge
d\overline{z}}{(1 + |z|^2)^2} \
\hskip 5mm = \ {\rm Kahler \ form \ on} \ \C
P^1.
\end{equation}

\subsubsection{$\cc_1$ of the determinant}\label{c1 det}

More generally, determinant bundles arise in geometric analysis,
in the representation theory of loop groups, and in the
construction of conformal field theories. In a general sense, they
facilitate the construction of projective representations from the
bordism category to categories of graded rings. The basic
invariant of a determinant bundle which one aims to compute is its
Chern class.


\subsubsection{Example: closed surfaces}\label{GRR}

A well known instance of that is for a family of compact
boundaryless surfaces $\{\Sigma_y\, |\, y\in Y\}$ parametrized by
a smooth manifold $Y$. Let $M = \bigcup_{y\in Y} \Sigma_y$ and
$\pi : M \to Y$ the projection map. Let $T_y$ be the tangent
bundle to $\Sigma_y$, and $T:= T(M/Y) =  \bigcup_{y\in Y} T_y \too
M$ the tangent bundle along the fibres. The index bundle $\Ind
\ol{\pd}_{(m)}$ of the family of D-bar operators $\ol{\pd}_{(m)} =
\{\ol{\pd}_y\ | \ y\in Y\}$ acting on sections of $T^{\otimes m}$
is the element $f_!(T^{\otimes m})$ of $K(Y)$, and the
Grothendieck-Riemann-Roch theorem says
\begin{equation*}
\ch(f_!(T^{\otimes m})) = f_*\left(\ch(T^{\otimes m})\, {\rm Todd}(T) \right),
\end{equation*}
where $f_* : H^i(M) \too H^{i-2}(Y)$ is integration over the fibres. That is,  with $\xi = c_1(T)$
\begin{eqnarray*}
\ch(\Ind \ol{\pd}_{(m)})  =  f_*\left(e^{m\xi} \cdot \frac{\xi}{1-e^{-\xi}} \right)
 =  f_*\left( 1+ (m + \frac{1}{2})\,\xi + \frac{1}{2}(m^2 + m +\frac{1}{6})\,\xi^2 + \ldots \right).
\end{eqnarray*}
Hence $\cc_1$ of the determinant line bundle $\Det\, \ol{\pd}_{(m)}$ is
\begin{equation}\label{GRR 2}
\cc_1 (\Det\, \ol{\pd}_{(m)} ) = \frac{1}{12}\,(6m^2 + 6m + 1)\,f_*(\xi^2) \ \in H^2(Y).
\end{equation}

\subsubsection{Quillen on the curvature formula}

More refined formulae may be sought at the level of smooth
invariants. The fundamental result in this direction was obtained
by Quillen in 1984 in a very beautiful four page article \cite{Qu}
in which the zeta function  regularized curvature of the
determinant line bundle $\Det \Ds_{\Sigma}$ of a family of Cauchy
Riemann operators $\Ds_{\Sigma} = \{D: \Omega(\Sigma,E) \to
\Omega^{0,1}(\Sigma,E)\}$ acting on sections of a complex vector
bundle $E$ over a closed Riemann surface $\Sigma$ was computed to
be
\begin{equation}\label{curv=kahler surface}
F_{\,\zeta}(\Ds_{\Sigma}) \ = \  {\rm  Kahler \ form \ on} \ Y
\end{equation}
where in this case $Y = \Omega^{0,1}(\Sigma,\End E)$.

\subsubsection{Bismut on Quillen}

Following Quillen's idea of constructing a superconnection on the
index bundle \cite{Qu2}, Bismut \cite{B} proved in a tour de force
a local index theorem for a general family $\Ds$  of Dirac-type
operators associated to a geometric fibration $\pi: M\too Y$ with
fibre a compact boundaryless manifold and, furthermore,  with
$F_{\,\zeta}(\Ds)\in\Omega^2(Y)$ the curvature of the
$\z$-connection on the determinant line bundle $\Det \Ds$,
extended \eqref{curv=kahler surface} to
\begin{equation}\label{curv=A hat ch}
F_{\,\zeta}(\Ds) \ = \   \textsf{ind}_{[2]},
\end{equation}
where $ \textsf{ind}\in\Omega^*(Y)$  is the family index density,
equal to $\int_{M/Y} \widehat{A}(M/Y)\, \ch(V)$ in the case of a
family of twisted Dirac operators, and the subscript indicates the
2-form component \cite{BF}.

It is worth emphasizing here the geometric naturality of the
formulae; in each of the above cases, including the example of
\secref{example complex projective space}, the $\zeta$-curvature
hits the index form `on the nose' --- any other connection will
have curvature differing from this by an exact 2-form.

\subsubsection{Melrose and Piazza on Bismut}\label{mp on b}

That naturality persists to the  analysis of families of APS
boundary problems $\Ds_{\Ps}$ for which the fibre of $\pi: M\too
Y$ is a compact manifold with boundary and $\pd M \neq \emptyset$,
and $\Ps = \{P_y\}$ is a smooth family of $\pdo$ projections on
the space of boundary sections which is pointwise (w.r.t. $Y$)
commensurable with the APS projection.

The principal contribution in this direction is the Chern
character formula of Melrose-Piazza \cite{MePi} proved using
$b$-calculus and generalizing Bismut-Cheeger \cite{BiCh}. From
this Piazza \cite{P} inferred the $b$ zeta-curvature function
formula on the $b$ determinant bundle $\Det^b(\Ds_\Ps)$ to be
\begin{equation*}
F^b_{\,\zeta}(\Ds_{\Ps}) \ = \ \textsf{ind}_{[2]}  \ + \ \left[\,\wt{\eta}_{\Ps}\,\right]_{(2)} ,
\end{equation*}
where $\wt{\eta}_{\,\Ps} := \pi^{-1/2}\int^{\oo}_0
\Tr(\dot{\Bs}_{t}\,e^{-\Bs_{t}^{2}}) \ dt$ is an eta-form of a
$t$-rescaled superconnection $\Bs_t = \Bs_t(\Ps)$ twisted by $\Ps$
for the family of Dirac  operators on the boundary $\pd M$.

\subsubsection{A direct computation}\label{here}

On the other hand, $\Ds_{\Ps}$ is already a smooth family of
Dirac-Fredholm operators and it is natural to seek a direct
computation of the $\zeta$-curvature formula for the determinant
line bundle $\Det\Ds_\Ps$, along the lines of  example of
\secref{example complex projective space}, without use of
$b$-calculus or other completions. It turns out, indeed, that
there is a canonical $\zeta$-function connection on $\Det\Ds_\Ps$
and one has:

\begin{thm}\label{thm intro}
Let $F_\z(\Ds_\Ps)$  be the curvature 2-form  of the $\zeta$-connection on $\Det\Ds_\Ps$. Then
\begin{equation}\label{det P curv intro}
F_\z(\Ds_\Ps) \ = \ F_\z(\Ds_{\Ps(\Ds)}) \ + \ \rkw  \hskip 10mm
in \ \ \Omega^2(Y)
\end{equation}
with  $\rkw$ the 2-form component of a relative $\eta$-form
depending only on boundary data; the fibration of closed boundary
manifolds and on ${\rm ran}(\Ps) = \Ww$ and on $\ran(\Ps(\Ds)) =
\Kk$. Here, $\Ps(\Ds)$ is the family of Calder\'{o}n projections
defined by $\Ds$, equal at $y\in Y$ to the projection onto the
(infinite dimensional) subspace equal to the restriction of  $\Ker
D_y$ to the boundary. The determinant bundle $\Det\Ds_{\Ps(\Ds)}$
is trivial. Its $\z$-curvature is canonically exact; there is a
preferred 1-form $\bet_\z(\Ds)\in \Omega^1(Y)$ such that
\begin{equation}\label{F calderon exact}
F_\z(\Ds_{\Ps(\Ds)}) = d\bet_\z(\Ds).
\end{equation}
\end{thm}

The definition of $\rkw$, which is simple and completely
canonical, and why it is a `relative eta  form', is given in
\secref{connection on sp}. The formula \eqref{det P curv intro} is
extremely `clean', in so far as it is the simplest relation that
might exist between $F_\z(\Ds_\Ps)$ and $\rkw$, both of which
represent $\cc_1(\Det\Ds_\Ps)$. It extends to geometric families
of boundary problems the principle of `reduction to the boundary'
present in the analysis of Grubb and Seeley \cite{GS1}, \cite{Gr2}
and Bruening and Lesch \cite{BL} of resolvent and zeta traces of
pseudodifferential boundary problems, also in Booss-Wojciechowki
\cite{BoWo}, and in the zeta determinant formulae in joint work
with Krzysztof Wojciechowski \cite{ScWo} and in \cite{Sc1}.



\subsubsection{Example: surfaces}\label{Gr and Diff S1 and CFT}

%

For a real compact surface $\Sigma$ with boundary $S^1$ our
conclusions  generalize the example of \secref{example complex
projective space} (and \secref{GRR}) as follows. A choice of
conformal structure $\tau\in {\rm {\rm Conf}}(\Sigma)$ turns
$\Sigma$ into a Riemann surface with a D-bar operator
$\ol{\pd}_\tau : \Omega^0(\Sigma) \too  \Omega^{0,1}(\Sigma)$.
Since  $P(\ol{\pd}_\tau)$ differs from the APS projection
$\Pi_\geq$  by only a smoothing operator \cite{Sc3} a suitable
parameter space of well-posed boundary conditions is the smooth
Grassmannian $\Gr$ of pseudodifferential operator ($\pdo$)
projections $P$ with $P - P(\ol{\pd}_\tau)$ smoothing. We obtain
in this way the family of APS boundary problems
$$\ol{\pd}_P := (\ol{\pd}_\tau)_P \ : \ {\rm dom} (\ol{\pd}_P) =
\Ker (P\circ \gam) \too  \Omega^{0,1}(\Sigma)$$
parametrized by $P\in\Gr$. In this case
$F_\z(\ol{\pd}_{P(\ol{\pd})})=0$ and \eqref{det P curv intro} is
\begin{equation}\label{curv=kahler Gr}
F_{\,\zeta}(\ol{\pd}_\Ps) \ = \ \Tr(PdPdP) \ = \  {\rm  Kahler \ form \ on} \ \Gr.
\end{equation}

The restriction of $\Det \ol{\pd}_\Ps$  to the loop group via the
embedding ${\rm LG}\hookrightarrow \Gr$  based at
$P(\ol{\pd}_\tau)$  is  the central extension  of ${\rm LG}$
(Segal \cite{Se}), while $F_{\,\zeta}(\ol{\pd}_\Ps)_{|{\rm LG}}$
is the 2-cocycle of the extension.
On the other hand, one may consider the opposite situation of  the
family of D-bar operators on $\Sigma$
$$\textsf{$\ol{\pd}$}_{\mbox{{\tiny$\Pi_\geq$,m}}} = \{(\ol{\pd}_\tau)_{\mbox{{\tiny$\Pi_\geq$,m}}} \ | \ \tau \in {\rm Conf}(\Sigma)\}$$
parametrized by $Y = {\rm Conf}(\Sigma)$ and with fixed boundary condition $\Pi_\geq$
acting on sections of $T^{\otimes m}\Sigma$. $\Det \textsf{$\ol{\pd}$}_{\mbox{{\tiny$\Pi_\geq$,m}}}$ pushes-down to  the moduli space $\Mm(\Sigma) = {\rm Conf}(\Sigma)/\Diff(\Sigma,\pd \Sigma)$ by the group of diffeomorphisms of $\Sigma$ equal to the identity on the boundary. In particular,  for the unit disc $D$ then $\Mm(D) = \Diff^+ S^1/{\rm PSU_{1,1}}$. By the functoriality of our constructions and the computations of \cite{Mo} we obtain that the $\z$-curvature of the determinant line bundle over $\Diff^+ S^1/{\rm PSU_{1,1}}$ is
\begin{equation*}
F_{\,\zeta}(\textsf{$\ol{\pd}$}_{\mbox{{\tiny$\Pi_\geq$,m}}}) \ =  \ F_\z(\textsf{$\ol{\pd}$}_{P(\ol{\pd}_m)}) \ + \  \frac{1}{12}\,(6m^2 + 6m + 1)\,\pi_*({\rm gv}) - \frac{1}{12}\,e,
\end{equation*}
where $\pi_*({\rm gv})$ is integration over the fibre of a Godbillon-Vey form, $e$ an Euler form \cite{Mo}, and $P(\ol{\pd}_m)$ the family of Calder\'{o}n boundary conditions.

\section{Fibrations of Manifolds}

Let $\pi : M \stackrel{X}{\too} Y$ be a smooth fibration of
manifolds with fibre diffeomorphic to a compact connected manifold
$X$ of dimension $n$ with boundary $\pd X\neq\emptyset$.
The total space $M$ is itself a manifold with boundary $\pd M$ and there is a boundary fibration $\sbp: \pd M \stackrel{\sbx}{\too}  Y$ of closed manifolds of dimension $n-1$. For example, for a fibration of surfaces over $Y=S^1$ then  $\pd M$ is a disjoint union of 2-tori fibred by the circle.

We assume there exists a collar neighbourhood $\Uu\< M$ of $\pd M$ with a diffeomorphism
\begin{equation}\label{collar M}
\Uu \cong [0,1) \times \pd M,
\end{equation}
corresponding fibrewise to a collar neighbourhood $[0,1) \times \pd X_y$ of each fibre $X_y :=\pi\ii(y) $.

\subsection{Bundles over fibrations}\label{bundles over fibrations}

A smooth family of vector bundles associated to $\pi : M \stackrel{X}{\too} Y$  is  defined to be a finite-rank $\Ci$ vector bundle
$E\too M.$
Formally, we may then consider the infinite-dimensional bundle $\Hh(E) \too Y$
whose fibre at $y\in Y$ is the space
$ \Hh_y(E) := \Gam(X_y,E_{|X_y})$
of $\Ci$ sections of $E$ over $X_y$. Concretely,  a section of $\Hh(E)$ is defined to be a  section of $E$ over $M$,
\begin{equation}\label{sections H(E)}
\Gam(Y,\Hh(E)) :\,=  \Gam(M,E).
\end{equation}
Thus, in practise one works with the right-side of \eqref{sections H(E)}, as indicated below.

$\Gam(Y,\Hh(E))$ is then a $\Ci(Y)$-module via
\begin{equation}\label{module}
\Ci(Y)\times\Gam(Y,\Hh(E)) \too \Gam(Y,\Hh(E)), \hskip 5mm (f,s) \mto f\cdot s :\,= \pi^*(f) s,
\end{equation}
that is, $f\cdot s(m) = f(\pi(m))\,s(m).$

The restriction map to boundary sections
\begin{equation}\label{restrict sections}
\gam: \Gam(Y,\Hh(E)) \too  \Gam(Y,\Hh(\Eb))
\end{equation}
is defined by the restriction map to the boundary on the total space
\begin{equation}\label{restrict sections 2}
\gam: \Gam(M,E) \too  \Gam(\pd M,\Eb)
\end{equation}
with $\Eb =  \cup_{m\in \pd M} E_m$ the bundle $E$ along $\pd M$.
Relative to  \eqref{collar M}
\begin{equation}\label{E on Uu}
E_{|\Uu} = \gam^*(\Eb)
\end{equation}
and
$
\Gam(\Uu,E) \cong \Ci([\,0,1))\otimes \Gam(M,\Eb).
$
Here, ${\rm rank}(\Eb) = {\rm rank}(E)$, so, for example, $TM_{{\tiny \partial M}}$  is not the same thing as $T(\pd M)$, whose sections are vector fields along the boundary, while a section of $TM_{{\tiny \partial M}}$ includes vector fields which point out of the boundary; one has $TM_{{\tiny \partial M}}\cong \R \oplus T(\pd M)$.

The vertical tangent bundle   $T(M/Y)$ (resp.  $T(\pd M/Y)$) is the subbundle of $TM$ (resp. $T\pd M$) whose fibre at $m\in M$ (resp. $m\in \pd M$) is the tangent space to the fibre $X_{\pi(m)}$ (resp. $\pd X_{\pi(m)}$).
 $\pi^*( TY)$  is the pull-back subbundle from the base.
Likewise, there is the dual bundle $T^* M$ with subbundle $T^*(M/Y)$, whose sections are vertical forms along $M$, and  $ \pi^*(\wedge T^* Y)$. More generally, the de-Rham algebra on $Y$ with values in $\Hh(E)$ is the direct sum of the
\begin{equation}\label{deR}
\Aa^k(Y,\Hh(E))  = \G(M,\pi^*(\wedge^k T^* Y)\otimes
E\otimes|\wedge_{\pi}|^{1/2}).
\end{equation}
The line bundle of vertical densities $|\wedge_{\pi}|$ is included to facilitate integration along the fibre.

\subsection{Connections}\label{connections}

A connection (or covariant derivative) on $\Hh(E)$ is specified by  a fibration `connection'  on $M$
\begin{equation}\label{M connection}
TM \cong T(M/Y) \oplus T_H M,
\end{equation}
and a vector bundle connection on $E$
\begin{equation}\label{e connection}
\wt{\nabla}:  \Gam(M,E) \to \Gam(M,E\otimes T^*M),
\end{equation}
which are compatible with the induced boundary connections.

The fibration connection is a complementary subbundle to $T(M/Y)$, specifying an isomorphism
$\pi^*(TY) \cong T_H M$ and hence a lift of vector fields from the base to horizontal vector fields on $M$
\begin{equation}\label{lift}
\G(Y,TY) \ \stackrel{\cong}{\too} \ \G(M, T_H M), \ \ \ \ \xi\mtoo \xi_H.
\end{equation}
A connection
\begin{equation}\label{nabla1}
\Nm: \Aa^0(Y,\Hh(E)) \too  \Aa^1(Y,\Hh(E))
\end{equation}
is then defined by
\begin{equation}\label{nabla2}
\Nm_\xi s = \wt{\nabla}_{\xi_H}s, \hskip 10mm s\in    \Gam(M,E), \ \xi\in\Ci(Y,TY).
\end{equation}

Compatibility with the boundary means, first, that  in the collar  $\Uu$
\begin{equation*}\label{connection in Uu}
\wt{\nabla}_{|\Uu} = \gam^*\tNb = \pd_u \,du + \tNb
\end{equation*}
where $u\in[\,0,1)$ is the normal coordinate to $\pd M$ and
$\tNb:  \Gam(\pd M,\Eb) \to \Gam(\pd M,\Eb\otimes T^*M)$
is the induced connection on  $\Eb$, defining $\Nb: \Aa^0(Y,\Hh(\Eb)) \too  \Aa^1(Y,\Hh(\Eb))$ by
\begin{equation}\label{Nb}
\Nb_\xi s = \tNb_{\xi_H}s, \hskip 10mm s\in \Gam(\pd M,\Eb).
\end{equation}
Secondly, that with respect to the boundary splitting
\begin{equation*}\label{N connection}
T(\pd M) \cong T(\pd M/Y) \oplus T_H \pd M
\end{equation*}
induced by
\begin{equation}\label{TUu}
T\Uu \cong \R \oplus T\pd M
\end{equation}
and the  splitting \eqref{M connection},
one has for $\xi\in\Ci(Y,TY)$ that
$$\left(\xi_H\right)_{|\Uu} \in\Ci(\Uu,T_H(\pd M)),$$
that is,
\begin{equation*}\label{assumption}
du(\xi_H) = 0,
\end{equation*}
where $du$ is extended from $\Uu$ to $M$ by zero.
One then has from \eqref{connection in Uu}
\begin{lem}
\begin{equation}\label{exact}
\gam\circ\wt{\nabla}_{\xi_H} = \wt{\nabla}^{\pd M}_{\xi_H}\circ\gam, \hskip 10mm \xi\in\Ci(Y,TY),
\end{equation}
as maps $\G(M,E) \too  \G(\pd M, \Eb)$.
\end{lem}

\vskip 1mm

The curvature of the connection \eqref{nabla1} evaluated on $\xi,\eta\in\Ci(Y,TY)$
\begin{equation}\label{curv nabla}
R(\xi,\eta) \in \G(Y, \,{\rm End}(\Hh(E))\,)
\end{equation}
is the smooth family of first-order  differential operators (as in \cite{B} Prop(1.11))
$$R(\xi,\eta) : = \wt{\nabla}_{\xi_H} \wt{\nabla}_{\eta_H} -\wt{\nabla}_{\eta_H}\wt{\nabla}_{\xi_H} - \wt{\nabla}_{[\xi,\eta]_H} = \wt{R}(\xi,\eta) +  \wt{\nabla}_{[\xi,\eta]_H- [\xi_H,\eta_H]}$$
where $\wt{R}(\xi,\eta)\in\G(M,\End E)$ is the curvature of $\wt{\nabla}$.
The above compatibility assumptions
state that $\gam_*(\xi_H),\gam_*(\eta_H)\in\Ci(\pd M, T_H(\pd M))$ and
\begin{equation}\label{curv nabla compatible}
 R(\xi_H,\eta_H)\circ \gam  = \Rb(\gam_*(\xi_H),\gam_*(\eta_H))\ \in \G(Y,\End(\Hh(\Eb))),
\end{equation}
where $\Rb(\alp,\bet)$ is  the curvature of \eqref{Nb}.

\subsubsection{Example: spin connection}\label{spinor}

For our purposes, here, it is not necessary to specify which particular connection on $E$ is being used, as the constructions are functorial. However, to compute the local index form curvature for a fibration of compact Riemannian spin manifolds (with or without boundary) then $\tNm$ must be the Bismut connection \cite{B}, \cite{BGV} and $E$ a twisted vertical spinor bundle. Then $T(M/Y)$ is oriented and spin, while  a metric $g$ on $TM$ in $\Uu$ is assumed to be the pull-back of a  metric $g^{\pd M}$  on $T\pd M$, so that $g^M_{|\Uu}  = du^2 + g^{\pd M}.$
If the connection on any twisting bundle is also of product type in the collar, then the situation of \secref{connections} holds, and the Bismut connection follows \cite{B}, \cite{BGV}.

\section{Families of Pseudodifferential Operators}

A smooth family of  $\pdos$ of constant order $\mu$ associated to
a fibration $\pi: N \stackrel{X}{\too} Y$ of compact boundaryless
manifolds, with $\dim(X)=n$,  with vector bundle $E\to N$ means a
classical $\pdo$ $$\As :\G(N,E^+)\too\G(N,E^-)$$ with Schwartz
kernel $k_\As\in\Dd^{'}(N\times_{\pi}N, E\boxtimes E)$ a vertical
distribution, where the fibre product $N\times_{\pi}N$ consists of
pairs $(x,x^{'})\in N\times N$ which lie in the same fibre, i.e.
$\pi(x) = \pi(x^{'})$,  such that in any local  trivialization
$k_\As$ is an oscillatory integral with vertical symbol $\as\in
S^{\nu}_{{\rm vert}}(N/Y)$ of order $\nu$. Here, $\xi$ is {\em
restricted} to the vertical momentum space, along the fibre. We
refer to $\As$ as a vertical $\pdo$ associated to the fibration
and denote this subalgebra of $\pdos$ on $N$ by
$$\G(Y,\Psi^\nu(E^+,E^-)) = \Psi^\nu_{{\rm vert}}(N,E^+,E^-).$$

In a similar way, for a fibration $\pi: M \too Y$ of compact
manifolds with boundary the pseudodifferential boundary operator
($\pdbo$) calculus as developed by Grubb \cite{Gr2}, generalizing
the Boutet de Monvel algebra, may be applied to define a vertical
calculus of operators with oscillatory integral kernels along the
fibres comprising trace operators from interior to boundary
sections, vertical Poisson operators taking  sections over the
boundary $\pd M$ into the interior, and restricted $\pdo$ and
singular Green's operators over the interior of $M$. This vertical
$\pdbo$ algebra is denoted
$$\G(Y,\Psi_\flat (E^+,E^-)) = \Psi_{{\rm vert},\flat}(M,E^+,E^-).$$

The algebras $\As\in \G(Y,\Psi^\nu(E^+,E^-))$ (see \cite{Sc2}) and
$\Psi_{{\rm vert},\flat}(M,E^+,E^-)$ of generalized $\pdos$ are
described in more detail in the Appendix.

For a local trivialization of the fibration and of $E$ one may
locally identify a vertical  $\pdo$ $\As$ with a single $\pdo$ (or
$\pdbo$) $A_y$ acting on a fixed space and depending on a local
parameter $y$ in $Y$.

\subsection{Families of Dirac-type operators}

Let $\Ds$ be a family of Dirac-type operators associated to the
fibration $\pi:M\to Y$ of compact manifolds with boundary with
vector bundles $E^\pm\to M$, such that in $\Uu$
\begin{equation}\label{D collar}
\Ds_{|\Uu} = \Upsilon\left(\frac{\pd}{\pd x_n} + \Dsb \right),
\end{equation}
where $\Dsb\in\Psi_{{\rm vert}}(\pd M, \Eb)$ a family of
Dirac-type operators associated to the boundary fibration of
closed manifolds, and $\Upsilon\in\G(\pd M, \End(\Eb))$ is a
bundle isomorphism.

\subsubsection{Vertical Poisson and Calder\'{o}n operators}

Let $\widehat{M} = M \cup_{\mbox{{\tiny $\pd M$}}}(-M)\to Y$ be
the fibration of compact boundaryless manifolds with fibre the
double manifold $\widehat{X}_y = X_y \cup_{\pd X_y} (-X_y)$. With
the product structure \eqref{collar M}, $\Ds$ extends by the proof
for a single operator, as in \cite{BoWo} Chap.9, to an {\em
invertible} vertical first-order differential operator
$\widehat{\Ds} \ \in \ \Psi^1(\wh{M},\wh{E}^+,\wh{E}^-),$ where
$\wh{E}^\pm_{|M} = E^\pm$ and $r^+\wh{\Ds}e^+ = \Ds$. As indicated
in  Appendix (A.1) and accounted for in detail in \cite{Sc2},
there is therefore a  smooth family of resolvent $\pdos$ of order
$-1$ $ \wh{\Ds}\ii \ \in \  \Psi_{{\rm
vert}}^{-1}(\wh{M},\wh{E}^-, \wh{E}^+). $ Define
\begin{equation*}
\Ds\ii_+ := r^+\wh{\Ds}\ii e^+ \ \in \  \Psi^{-1}_{{\rm
vert},\,\flat}(M,E^-, E^+).
\end{equation*}
Since $\widehat{\Ds}\widehat{\Ds}\ii = \Is$  on
$\G(\widehat{M},\widehat{\Ee})$, with $\Is$ the vertical identity
operator, and since $\Ds$ is local
\begin{equation}\label{DDii}
\Ds\Ds\ii_+  = \Is \ \ \ \ \ {\rm on} \ \ \G(M,E^-).
\end{equation}
Thus there is a short exact sequence
$0 \too \Ker(\Ds) \too \G(M,E^+) \stackrel{\Ds}{\too}
\G(M,E^-)\too 0$, where
\begin{equation}\label{ker Ds}
\Ker(\Ds) = \{s\in \G(M,E^+) \ | \ \Ds s = 0 \ \
{\rm in} \ \ M\backslash \pd M\}.
\end{equation}
On the other hand,  $\Ds\ii_+$ is not a left-inverse but (by an obvious modification of \cite{S1}, \cite{S2}, \cite{BoWo} \S 12)
\begin{equation}\label{Dinverse D}
\Ds\ii_+\Ds  = \Is - \Ks \gam\ \ \ \ \ {\rm on} \ \ \G(M,E^+),
\end{equation}
where $\gam$ is the restriction operator \eqref{restrict sections 2} and
the  {\it vertical Poisson operator associated to  $\Ds$} is
\begin{equation}\label{Poisson}
\Ks = \Ds\ii_+ \gam^* \Upsilon,
\end{equation}
with $\gam$ as in \eqref{poisson operators}. Composing with boundary restriction defines the {\it vertical Calder\'{o}n projection} (\cite{Ca}, \cite{S1}, \cite{S2}, \cite{BoWo})
\begin{equation}\label{calderon}
P(\Ds) := \gam\circ \Ks \ \  \in \ \G(Y, \Psi_{{\rm vert}}^0(\Eb)):\,= \ \Psi_{{\rm vert}}^0(\pd M, \Eb)
\end{equation}
with range the space of vertical Cauchy data
\begin{equation}\label{ran calderon}
\ran(P(\Ds)) = \gam \Ker(\Ds) = \{f\in \G(\pd M, \Eb) \ | \ f=
\gam s, \ s \in \Ker(\Ds)\}.
\end{equation}
This may be formally characterized as the space of sections of the infinite-dimensional subbundle
$\Kk(\Ds)\< \Hh(\Eb)$
with fibre $K(D_y) = \gam\Ker(D_y)$ at $y\in Y$ (and, likewise, $\Ker\Ds$ as the space of sections of the formal subbundle of $\Hh(E^+)$ with fibre $\Ker D_y$). However, as with  $\Hh(\Eb)$ in \secref{bundles over fibrations}, \,concretely one only  works with the space of sections of $\Kk(\Ds)$
\begin{equation}\label{ran calderon k}
 \G(Y,\Kk(\Ds))  :\,= \{f\in \G(\pd M, \Eb) \ | \ f=
\gam s, \ s \in \Ker(\Ds)\} =\ran(P(\Ds)).
\end{equation}
(Note, on the other hand, $\Kk(\Ds)$ is not the space of sections of a subbundle of $\Eb$.)

By the fibrewise Unique Continuation property, restriction $\gam: \Ker \Ds \too \G(Y,\Kk(\Ds))$ defines a canonical isomorphism with right-inverse
\begin{equation}\label{KerDs=Kk(Ds)}
\Ks : \G(Y,\Kk(\Ds)) \stackrel{\cong}{\too} \Ker(\Ds).
\end{equation}

\subsection{Well-posed boundary problems for D}

The vertical Calder\'{o}n projection \eqref{calderon} provides the
reference $\pdo$ on boundary sections with respect to which is
defined any vertical well-posed boundary condition for  $\Ds$.

\subsubsection{Smooth families of boundary $\pdo$ projections}\label{boundary projs}

We consider smooth families of $\pdos$ on $\G(Y,\Hh(\Eb))$ which
are perturbations of the Calder\'{o}n projection of the form
\begin{equation}\label{Ps}
  \Ps = P(\Ds) + \verb"S" \ \ \in \Psi^0_{{\rm vert}}(\pd M, \Eb),
\end{equation}
where $$\verb"S" \in \Psi^{-\oo}_{{\rm vert}}(\pd M, \Eb)$$ is a
vertical smoothing operator (smooth family of smoothing
operators), cf. Appendix. From Birman-Solomyak \cite{BiSo}, Seeley
\cite{S2} (see also \cite{BoWo}) \eqref{Ps} may be replaced by the
projection  onto $\ran(\Ps)$ to define an equivalent boundary
problem. So we may assume $\Ps^{\,2} = \Ps$ and $\Ps^* = \Ps,$
where the adjoint is with respect to the Sobolev completions and
vertical inner-product defined by metric on $\Eb$ and the choice
of vertical density  $d_{\mbox{{\tiny $\pd$ M/Y}}}x^{'} \in\G(\pd
M,  |\wedge^{n-1} T^*(\pd M/Y)|)$.

The family APS projection $\Pi_{>} = \{\Pi^{y}_{>} \ | \ y\in Y\}$
is only smooth in $y$ when $\dim\Ker(\Dsb)_y$ is constant
\cite{BiCh}. Nevertheless, we refer to \eqref{Ps} as a vertical
$\pdo$ of APS-type.

The choice of $\Ps$ in \eqref{Ps} distinguishes the subspace of
the space of  boundary sections
\begin{equation}\label{Ww}
\G(Y,\Ww) := \ran(\Ps) = \{\Ps f \ | \ f\in \G(\pd M,\Eb) \} \ \< \ \Gam(\pd M, \Eb) := \G(Y,\Hh(\Eb)).
\end{equation}
Here, $\Ww$  is the formal infinite-rank subbundle of $\Hh(\Eb)$ with fibre $W_y =\ran \Ps_y\<\G(\pd X_y, (\Eb)_y)$, whose local bundle structure follows from the invertibility of the operators $P_{y^{'}}P_y : W_y \too W_{y^{'}}$ for $y^{'}$ near $y$. Analytically, though, just as with $\Kk(\Ds)$, one works in practise with \eqref{Ww}.

Given any two choices $\Ps, \Ps^{'}$ of the form \eqref{Ps} one
has the smooth family of Fredholm operators
\begin{equation}\label{P2oP1}
\Ps^{'} \circ \Ps : \G(Y,\Ww) \too  \G(Y,\Ww^{'})
\end{equation}
where $\G(\pd M,\Ww) :\, = \ran(\Ps)$. We may write this as a
section of the formal bundle $\Hom(\Ww,\Ww^{'})$ in so far as we
declare the sections of the latter to precisely be the subspace of
$\Psi_{{\rm vert}}(\pd M, \Eb)$
\begin{equation}\label{Hom(W1,W2)}
\G(Y, \Hom(\Ww,\Ww^{'})) :\,= \{\Ps^{'} \circ\As\circ \Ps \ | \
\As\in \Psi_{{\rm vert}}(\pd M, \Eb)\}.
\end{equation}
Note here that
$$\Ps^{'} \circ \Ps \in \Psi^0_{{\rm vert}}(\pd M, \Eb)$$
is a smooth vertical $\pdo$ on boundary sections. The reference to it as a `smooth family of Fredholm operators' means additionally that there is smooth vertical $\pdo$ on boundary sections
$$\Qs_{\mbox{{\tiny$\Ps,\Ps^{'}$  }}} \in \Psi^0_{{\rm vert}}(\pd M, \Eb)$$ such that
\begin{equation}\label{QP2oP1}
\Qs_{\mbox{{\tiny$\Ps,\Ps^{'}$  }}}\circ(\Ps^{'} \circ \Ps)  = \Ps
+ \Ps\Sss^{'}\Ps,   \hskip 10mm \Sss^{'}\in \Psi^{-\oo}_{{\rm
vert}}(\pd M, \Eb),
\end{equation}
and hence that $\Qs_{\mbox{{\tiny$\Ps,\Ps^{'}$  }}}$ is a
parametrix for \eqref{P2oP1}; that is, restricted to $\G(Y,\Ww)$
\eqref{QP2oP1} is
\begin{equation}\label{QP2oP1onWw1}
\left(\Qs_{\mbox{{\tiny$\Ps,\Ps^{'}$  }}}\circ(\Ps^{'} \circ
\Ps)\right)_{|\Ww} = \Is_{\mbox{{\tiny $\Ww$}}} + \Ps\Sss^{'}\Ps,
\end{equation}
where $\Is_{\mbox{{\tiny $\Ww$}}}$  denotes the identity on
$\G(Y,\Ww)$. Indeed, we may take, for example,
$\Qs_{\mbox{{\tiny$\Ps,\Ps^{'}$  }}} = \Ps \circ \Ps^{'}.$

\subsubsection{Vertical APS-type boundary problems}\label{Dp}

The choice of $\Ps$ in \eqref{Ps} additionally distinguishes the subspace of interior sections on the total space of the fibration (which is not itself the space of sections of some subbundle of $E^+$)
\begin{equation}\label{Hp}
\G(Y,\Hh_\Ps(E^+)) :\;= \Ker(\Ps \circ \gam) = \{s\in \G(M,E^+) \ | \ \Ps\gam s=0\} \ \< \ \G(M,E^+) :\;= \G(Y,\Hh(E^+)).
\end{equation}
We may consider the infinite-dimensional bundle $\Hh_\Ps(E) \too Y$ with fibre at $y\in Y$  the space of $\Ci$ sections of $E^+$ over $X_y$ which lie in $\Ker (P_y \circ \gam)$, related to $\Ww$ via the exact
sequence \
$0 \too \Hh_\Ps(E^+) \too \Hh(E^+) \stackrel{\Ps\circ \gam
}{\too} \Ww \too 0$. Concretely, however, one works in practise  with \eqref{Hp}.

A smooth family of APS-type boundary
problems is the restriction of $\Ds$ to the subspace \eqref{Hp}
\begin{equation}\label{DP}
\Ds_{\Ps}\,:=\,\Ds\,:\ \Ker(\Ps \circ \gam) = \G(Y,\Hh_\Ps(E^+))\too \Gam(M,E^-).
\end{equation}

\vskip  1mm  $\Ds_\Ps$ restricts over $X_y$ to $D_{P_y} := (D_y)_{P_y} :{\rm dom}(D_{P_y})
\too \G(X_y,E^-_y)$ in a local trivialization of the fibration of manifolds, an
APS boundary problem in the usual single operator sense.

The existence of the Poisson operator \eqref{Poisson} reduces the construction of
a vertical parametrix for $\Ds_{\Ps}$ to the construction of a parametrix for the
operator \eqref{P2oP1} on boundary sections
\begin{equation}\label{PoP(Ds)}
\sp := \Ps \circ P(\Ds) : \ \G(Y,\Kk(\Ds)) \too \G(Y,\Ww).
\end{equation}
Explicitly, Let $U\<Y$ be the open subset of points in $Y$ where
$\sp$ is invertible. That is,  relative to any local
trivialization of the geometric fibration $M\to Y$ and bundles at
$y\in Y$ the Fredholm family $S(P)$ parametrizes an operator
$S_y(P_y)=P_y\circ P(D_y) : K(D_y) \too \ran(P_y)$ in the usual
single operator sense; $y\in U$ if $S_y(P_y)$ is invertible. Over
$U$ we define
\begin{equation}\label{Poisson Dp}
\Ks(\Ps)_{|U} := \Ks \circ P(\Ds) \sp_{|U}\ii \Ps : \
\G(\pi_{\mbox{{\tiny $\pd$} }}\ii(U), \Eb) \too
\G(\pi_{\mbox{{\tiny $\pd$} }}\ii(U),E^+),
\end{equation}
where $\pi_{\mbox{{\tiny $\pd$} }} : \pd M \to Y$ is the boundary
fibration. Then  Green's theorem for the vertical densities along
the fibres locally refines \eqref{Dinverse D} to
\begin{equation}\label{DPmaster}
(\Ds_\Ps)_{|U}\ii \Ds = \Is_{|U} - \Ks(\Ps)_{|U}\gam : \
\G(\pi_{\mbox{{\tiny $\pd$} }}\ii(U),E^+) \too
\G(\pi_{\mbox{{\tiny $\pd$} }}\ii(U),E^+).
\end{equation}
Moreover, if $\Ds_{\Ps^{'}}$ is also invertible over $U$
\begin{equation}\label{relinv}
(\Ds_{\Ps})_{|U} \ii =
(\Ds_{\Ps})_{|U}\ii\Ds(\Ds_{\Ps^{'}})_{|U}\ii = \Ds_{\Ps^{'}}\ii -
\Ks(\Ps)_{|U}\Ps\gam \Ds_{\Ps^{'}} \ii : \ \G(\pi_{\mbox{{\tiny
$\pd$} }}\ii(U),E^-) \too \G(\pi_{\mbox{{\tiny $\pd$}
}}\ii(U),E^+),
\end{equation}
We note, globally on $M$, that:
\begin{prop}\label{relatively smooth}
With the above assumptions the relative inverse is a vertical
smoothing operator
\begin{equation}\label{relinv smooth}
(\Ds_{\Ps})_{|U}\ii - (\Ds_{\Ps^{'}})_{|U}\ii   \ \ \in \ \G(U,
\Psi^{-\oo}_{{\rm vert,\flat}}(E_{|\pi_{\mbox{{\tiny $\pd$}
}}\ii(U)})).
\end{equation}

More generally, for a general APS-type vertical $\pdo$ projection
$\Ps\in\Psi^0_{{\rm vert}}(\pd M, \Eb)$ a global parametrix for
the smooth family of APS-type boundary problems  $\Ds_\Ps:
\Ker(\Ps \circ \gam)\too \Gam(M,E^-)$ is given by
\begin{equation}\label{DPparametrix}
\Ds_+\ii - \Ks\, \Qs_{\mbox{{\tiny$\Ps,P(\Ds)$  }}} \g \Ds_+\ii \
\ \in \ \G(Y, \Psi^{-\oo}_{{\rm vert,\flat}}(E_{|M})),
\end{equation}
where $\Qs_{\mbox{{\tiny$\Ps,P(\Ds)$  }}}$ is any parametrix as in
\eqref{QP2oP1onWw1} for $\sp$, for example
$\Qs_{\mbox{{\tiny$\Ps,P(\Ds)$  }}} = P(\Ds)\circ\Ps$.
\end{prop}
{\bf Proof.} \ We have $\Ps= P(\Ds) + \verb"S"$, $\Ps^{'}= P(\Ds)
+ \verb"S"^{'}$ for vertical smoothing operators
$$\verb"S", \ \verb"S"^{'}\in \Psi^{-\oo}_{{\rm vert}}(\pd M, \Eb).$$
Hence
\begin{equation}\label{P1-P2 a}
\Ps - \Ps^{'} \in \Psi^{-\oo}_{{\rm vert}}(\pd M, \Eb)
\end{equation}
and
\begin{equation}\label{P1-P2}
\Ps(\Is-\Ps^{'}) = -\Ps\,\verb"S"^{'} \in \Psi^{-\oo}_{{\rm vert}}(\pd M, \Eb)
 \end{equation}
are  vertical smoothing operator operators. By \eqref{relinv}
\begin{equation}\label{relinv2}
(\Ds_{\Ps})_{|U}\ii - (\Ds_{\Ps^{'}})_{|U}\ii = - \Ks(\Ps)\Ps\gam
(\Ds_{\Ps^{'}})_{|U} \ii =
 - \Ks(\Ps)_{|U}\Ps(\Is-\Ps^{'})\gam(\Ds_{\Ps^{'}})_{|U} \ii \ \ \ \ {\rm over} \ \ M_{|U} = \pi\ii(U)
\end{equation}
which by \eqref{P1-P2} and the composition rules of the $\pdbo$ calculus (cf \secref{families pdbos}\,) is smoothing.

The assertion that \eqref{DPparametrix} is a parametrix is an obvious slight modification of the argument leading to \eqref{relinv}.
\begin{flushright}
$\Box$
\end{flushright}

\section{The Determinant Line Bundle}

From \propref{relatively smooth} the choice of $\Ps$ restricts $\Ds$ to a
family $\Ds_\Ps$ of Fredholm operators. It also has the consequence that the
kernels of the restricted operators no longer  define a vector
bundle (formally \eqref{ker Ds} does), rather they define a virtual bundle $\Ind \Ds_\Ps  \in  K(Y)$.
Likewise, from \secref{boundary projs},  $\sp : \G(Y,\Kk((\Ds)) \too \G(Y,\Ww)$ is a smooth Fredholm family defining an element $\Ind \sp  \in  K(Y)$. The determinant line bundles $\Det \Ds_\Ps$ and $\Det \sp$ are the top exterior powers of these elements, at least in K-theory. To make sense of them as smooth complex line bundles we use the following trivializations, with respect to which the zeta connection will be constructed.

\subsubsection{Determinant lines}

The determinant of a Fredholm operator $ T:H\to H^{'}$ exists
abstractly not as a number but as an element $\det T$ of a complex
line $\Det T$. A point of $\Det T$
is an equivalence class $[S,\lam]$ of pairs $(S,\lam)$, where
$S:H\to H^{'}$ differs from $T$ by a trace-class
operator and relative to the equivalence relation $(Sq,\lam) \sim
(S,\lam\det_F q)$ for $q:H \to H$ of Fredholm-determinant class.
Scalar multiplication on $\Det T$ is  $\mu . [S,\lam] = [S,\mu\lam].$
The determinant $\det T := [T,1]$ is non-zero if and only if $T$ is invertible,
and there is a canonical isomorphism
\begin{equation}\label{det line}
\Det T \cong \wedge^{{\rm max}}\Ker T^*
\otimes \wedge^{{\rm max}}\Cok T.
\end{equation}
For Fredholm operators $T_1, T_2 :H\to H^{'}$ with $T_i - T$ trace class and $T_2$ invertible
\begin{equation}\label{det ratio}
  \frac{{\rm det}\, T_1}{{\rm det} \,T_2} =  {\rm det}_F (
  T_1T_2\ii),
\end{equation}
where the quotient on the left side is taken in $\Det T$ and ${\rm det}_F $ on the right-side in $H^{'}$.

\subsubsection{The line bundle Det\,$\sp$}\label{det sp}

For each smooth family of smoothing operators $\s =\{\s_y\}\in \G(Y, \Psi^{-\oo}_{{\rm vert}}(\Eb)) = \Psi^{-\oo}_{{\rm vert}}(\pd M, \Eb)$ define
\begin{equation}\label{Pss}
\Pss = \Ps + \Ps\s\Ps\in \G(Y, \Psi^0_{{\rm vert}}(\Eb)) = \Psi^0_{{\rm vert}}(\pd M, \Eb)
\end{equation}
and the open subset of $Y$
\begin{equation}\label{Us}
U_{\s} := \{y\in Y \ | \ \Sss(\Pss)_y:= (P_y+ P_y\,\sigma_y \,P_y)\circ P(D_y) : K(D_y) \too \ran(P_y) \ \ {\rm invertible}\}.
\end{equation}
Over $U_\s$ one has the canonical trivialization
\begin{equation}\label{triv Us}
U_\s \ \too \Det \sp_{|U_\s} = \bigcup_{y\in U_\s} \Det\sp_y, \ \ \ \ y\mtoo \det\Sss(\Pss)_y : = [ (\Pss\circ P(\Ds))_y,1],
\end{equation}
where $\sp_y:=P_y\circ P(D_y):K(D_y)\to \ran(P_y)$.
Note that $\det \Sss(\Pss)_y\neq 0$, and that $\Sss(\Pss)_y - \sp_y$ is the restriction of a smoothing operator so that
\begin{equation}\label{spsigma in Det sp}
\det \Sss(\Pss)_y\ \in \ \Det\sp_y\bsh\{0\}.
\end{equation}
Over the intersection $U_{\s}\cap U_{\s^{'}}\neq\emptyset$ the transition function by \eqref{det ratio} is the function
\begin{equation}\label{transition Det sp}
U_{\s}\cap U_{\s^{'}} \ \too \ \C^*, \ \ \ y \mtoo {\rm det}_F \left(\Sss(\Pss)_y \circ \Sss(\Ps_{\mbox{{\tiny $\sigma'$}}})_y\ii \right),
\end{equation}
where the Fredholm determinant is taken on $\ran(P_y)$ and varies holomorphically with $y$.

\subsubsection{The line bundle Det\,$\Ds_\Ps$}\label{det Dp}

The bundle structure of $\Det\Ds_\Ps$ is defined by perturbing $D_{P_y}$ to an invertible operator. It is crucial for the construction of the $\z$ connection to do so by perturbing the $\pdo$ $P_y$, not $D_y$.

To do this we mediate the local trivializations of $\Det\Ds_\Ps$ through those of $\Det\sp$ in \secref{det sp}.

Precisely, the family of $\pdos$ $\Pss$ in \eqref{Pss} is of APS-type
\begin{equation}\label{aps type}
\Pss -  P(\Ds)  \ \ \in \Psi^{-\oo}_{{\rm vert}}(\pd M, \Eb),
\end{equation}
defining the vertical boundary problem $\Ds_{\,\Pss} : \G(Y,\Hh_{\Pss}(E^+)) \too \G(Y,\Hh(E^-))$. From \secref{Dp}
\begin{equation}\label{Us 2}
U_\sigma := \left\{y\in Y \ | \ (\Ds_{\,\Pss})_y :{\rm dom}((D_{\Pss})_y
\too \G(X_y,E^-_y) \ \ {\rm invertible}\right\},
\end{equation}
over which there is the local trivialization
\begin{equation}\label{triv Us 2}
U_\sigma \ \too \Det \Ds_{\Pss|U_s}, \ \ \ \ y\mtoo \det((\Ds_{\Pss})_y) = [(\Ds_{\Pss})_y, 1] \in \Det (\Ds_{\Pss})_y.
\end{equation}
The equivalence of \eqref{Us} and \eqref{Us 2} is the  identification for any APS-type  $\wt{P}\in \Psi^0_{{\rm vert}}(\pd M, \Eb)$ of the kernel of $(D_{\wt{P}})_y $ with that of $S(\wt{P})_y$ defined by the Poisson operator $\Ks_y$, and likewise of the cokernels. It follows that there is a {\em canonical} isomorphism
\begin{equation}\label{Det dp = Det sp}
\Det (D_{\wt{P}})_y \ \cong \ \Det S(\wt{P})_y \ \ \ {\rm with} \ \ \ \det (D_{\wt{P}})_y \longleftrightarrow \det S(\wt{P})_y.
\end{equation}
The local trivialization of $\Det \Ds_\Ps$ is then defined through the canonical isomorphisms of complex lines applied to \eqref{triv Us 2}
\begin{equation}\label{canonical isoms}
 \Det (D_{\Pss})_y  \stackrel{\eqref{Det dp = Det sp}}{\cong} \Det\Sss(\Pss)_y = \Det\sp_y \stackrel{\eqref{Det dp = Det sp}}{\cong} \Det (\Ds_{\Ps})_y,
\end{equation}
where the central equality is from \secref{det sp}. By construction the transition functions for  $\Det \Ds_\Ps$ are precisely \eqref{transition Det sp}; that is, as functions of $y\in U_{\s}\cap U_{\s^{'}}$
\begin{equation}\label{transition Det Dp}
\det (\Ds_{\Pss})_y = {\rm det}_F \left(\Sss(\Pss)_y \circ \Sss(\Ps_{\mbox{{\tiny $\sigma'$}}})_y\ii \right)\,\det (\Ds_{\,\Ps_{\mbox{{\tiny $\sigma'$}}}})_y\ \ \ \ \ {\rm in} \ \ \Det(\Ds_\Ps)_y.
\end{equation}

Thus the bundle structure of $\Det \Ds_\Ps$ is constructed using that of $\Det\sp$, as with all other spectral invariants of $\Ds_\Ps$ owing to the facts in \secref{Dp}.

With respect to  smooth families of boundary conditions  $\Ps,
\Ps^{'} \in \Psi^0_{{\rm vert}}(\pd M, \Eb)$
\begin{equation}\label{detbundle P}
\Det \Ds_{\Ps} \cong  \Det \Ds_{\Ps^{'}} \otimes \Det
(\Ps\circ\Ps^{'}),
\end{equation}
which may be viewed as a smooth version of the K-theory identity
\begin{equation}\label{indexbundle P}
\Ind \Ds_{\Ps} = \Ind \Ds_{\Ps^{'}} + \Ind (\Ps\circ\Ps^{'}).
\end{equation}
These are a consequence of the following general (useful) identifications.
\begin{thm}\label{mult in bundle}
Let $\As_1:\Hh^{'}\to \Hh^{''}$, $\As_2:\Hh\too \Hh^{'}$ be smooth (resp. continuous) families of Fredholm operators acting between Frech\'{e}t bundles over a compact manifold $Y$. Then there is a canonical isomorphism of $\Ci$ (resp $C^0$) line bundles
\begin{equation*}
\Det \As_1\As_2 \cong  \Det \As_1 \otimes
\Det \As_2
\end{equation*}
with $\det \As_1\As_2  \longleftrightarrow \det\As_1\otimes \det\As_2$. In $K(Y)$ one has
\begin{equation}\label{ind P}
\Ind \As_1\As_2 =  \Ind \As_1 \ + \ \Ind \As_2
\end{equation}
\end{thm}

For a proof of \thmref{mult in bundle} see \cite{Scott book}.

\section{Hermitian Structure}

The (Quillen) $\z$-metric on $\Det \Ds_\Ps$ is defined  over
$U_\s$ by evaluating it on the  non-vanishing section $\det
\Ds_{\Pss}$
\begin{equation}\label{z metric}
\|\det (\Ds_{\Pss})_y\|_\z^2 =  {\rm det}_{\zeta}(\Delta_{\Pss})_y,
\end{equation}
where the right-side is the $\z$-determinant of the  vertical
Laplacian boundary problem for an APS-type $\pdo$ $\Ps$
\begin{equation}\label{Delta P}
\Delta_\Ps = \Delta := D^*D \ : \ {\rm dom}(\Delta_\Ps) \to \G(M,E^-)
\end{equation}
with ${\rm dom}(\Delta_\Ps) = \{s\in\G(M,E^+) \ | \ \Ps\gam s=0, \
\Ps^\star\gam Ds=0 \}$ and $\Ps^\star :=\Upsilon (I-
P_y)\Upsilon^*$ the adjoint vertical boundary condition.

From \cite{Sc1} Thm(4.2) we know that
\begin{equation}\label{z metric patching}
\|\det (\Ds_{\Pss})_y\|_\z^2  = \frac{{\rm det}_F \left(\Sss(\Pss)_y^*\, \Sss(\Pss)_y\right)}{{\rm det}_F\left( \Sss(\Ps_{\mbox{{\tiny $\sigma'$}}})_y^*\,\Sss(\Ps_{\mbox{{\tiny $\sigma'$}}})_y\right)}\,\|\det (\Ds_{\,\Ps_{\mbox{{\tiny $\sigma'$}}}})_y\|_\z^2,
\end{equation}
which is the patching condition with respect to the transition functions \eqref{transition Det Dp} for \eqref{z metric} to define a global metric on the determinant line bundle $\Det \Ds_\Ps$.

\section{Connections on Det$\Ds_\Ps$}\label{connection on sp}

There are two natural  ways to put a connection on the determinant bundle $\Det \Ds_\Ps$. The first of these is associated to the boundary fibration and its curvature may be viewed as a relative $\eta$-form. The second, is the $\z$-function connection, the object of primary interest here.

\subsection{A connection on $\Det\sp$ }\label{connection on Det sp}

The first connection is defined on $\Det \sp$, which defines a connection on $\Det \Ds_\Ps$ via the isomorphism (by construction) between these line bundles.

The endomorphism bundle $\End(\,\Hh(\Eb)) $ whose sections are the boundary vertical $\pdos$
$$\G(Y,\End(\,\Hh(\Eb)))  :\,= \Psi^*_{{\rm vert}}(\pd M, \Eb)$$
has an induced connection (also denoted $\nablab$) from $\nablab$ on $\G(Y,\Hh(\Eb))$ in \eqref{Nb} by
$$\nablab_\xi\As := [\,\nablab_\xi,\,\As] \ \in \Psi^*_{{\rm vert}}(\pd M, \Eb),$$
where $\xi\in\Ci(Y,TY)$. That is,
\begin{equation}\label{nabla hom}
(\nablab_\xi \As) f  = \tNb_{\xi_H}(\As f) - \As(\tNb_{\xi_H} f), \hskip 10mm f\in \G(\pd M,\Eb).
\end{equation}
Let $P(\Ds)\in \Psi^0_{{\rm vert}}(\pd M, \Eb)$ be the Calder\'{o}n vertical $\pdo$ projection, and let  $\Ps\in \Psi^0_{{\rm vert}}(\pd M, \Eb)$ be any other vertical APS-type boundary condition \eqref{Ps}. Then there are  induced connections
$$\Nbw = \Ps\cdot\Nb\cdot\Ps, \hskip 10mm \Nbk = P(\Ds)\cdot\Nb\cdot P(\Ds)$$ defined on the Frech\'{e}t bundles $\Ww$ and $\Kk(\Ds)$, in the sense that
$$\Nbw_\xi : \  \G(Y,\Ww) :\,= \{\Ps s \ | \ s\in \G(\pd M, \Eb)\} \ \too \ \G(Y,\Ww)$$
with
$$\Nbw_\xi s = \Ps\,\tNb_{\xi_H}(\Ps s), \hskip 10mm s\in \G(Y,\Ww),$$
satisfies the Leibnitz rule, and likewise for $\Nbk$. We therefore have the induced connection $\Nbkw$ on the restricted hom-bundle $\Hom(\Kk(\Ds),\Ww)$, where, as in \eqref{Hom(W1,W2)},
\begin{equation}\label{hom w k}
\G(Y,\Hom(\Kk(\Ds),\Ww)) \ :\,= \ \{\Ps\circ\Cs\circ P(\Ds) \ | \ \Cs \in \Psi^*_{{\rm vert}}(\pd M, \Eb)\},
\end{equation}
defined by
\begin{equation}\label{nabla hom k w}
(\Nbkw_\xi \As) s = \Nbw_\xi(\As s) - \As(\Nbk_\xi s) , \hskip 10mm s\in \G(\pd M,\Eb), \ \As\in \G(Y,\Hom(\Kk,\Ww)).
\end{equation}
One then has a connection on $\Det \sp$ by setting over $U_\s$
\begin{equation}\label{nabla det k w}
\Nsp_{|U_\s} \det \Sss(\Pss) =  \omega^{\,\mbox{{\tiny $\Sss(\Ps_\s)$}}}\,\det \Sss(\Pss)
\end{equation}
where the locally defined 1-form in $\Omega^1(U_\s)$ is
\begin{equation}\label{omega det k w}
\omega^{\,\mbox{{\tiny $\Sss(\Ps_\s)$}}} \ = \  \Tr(\Sss(\Ps_\s)\ii\nabla_\xi^{\,\mbox{{\tiny ${\mathcal K},{\mathcal W}_\sigma$}}}\,\Sss(\Ps_\s)),
\end{equation}
with $\G(Y,\Ww_\sigma) = \ran(\Ps_\s)$. The trace on the right-side of \eqref{omega det k w} is the usual vertical trace (along the fibres, as recalled in the Appendix), by construction taken over $\G(Y,\Kk(\Ds))\<\,\G(\pd M, \Eb)$.

Notice, here, that $\As\ii\Nbkw_\xi \As$ will not be a {\em trace-class} family of $\pdos$ for a general invertible vertical $\pdo$ $\As\in \G(Y,\Hom(\Kk,\Ww))\<\Psi^*_{{\rm vert}}(\pd M, \Eb)$. That this is nevertheless the case when $\As = \Ss(\Ps)$, so that the right-side of \eqref{omega det k w} is well defined, is immediate from \eqref{aps type} and \eqref{nabla hom k w}.

The local 1-forms define a global connection with respect to \eqref{transition Det sp} by the identity in $\Omega^1(U_\s\cap U_{\sigma'})$
\begin{equation*}
d_\xi\, {\rm det}_F \left(\Sss(\Pss) \circ \Sss(\Ps_{\mbox{{\tiny $\sigma'$}}})\ii \right)\ = \  \Tr(\Sss(\Ps_\s)\ii\nabla_\xi^{\,\mbox{{\tiny ${\mathcal K},{\mathcal W}_\sigma$}}}\Sss(\Ps_\s)) - \Tr(\Sss(\Ps_{\mbox{{\tiny $\sigma'$}}})\ii\nabla_\xi^{\,\mbox{{\tiny ${\mathcal K},{\mathcal W}_{\sigma'}$}}}\Sss(\Ps_{\mbox{{\tiny $\sigma'$}}})).
\end{equation*}
which is a standard Fredholm determinant identity $d_\xi\, {\rm det}_F\Cs = \Tr(\Cs\ii\nabla_\xi\Cs)$ for a smooth family of Fredholm-determinant class operators $y\mto\Cs(y)$.

\subsubsection{Curvature of $\Nsp$ }

The curvature of the connection  $\Nsp$ on the complex line bundle $\Det \sp \too Y$ is the globally defined 2-form
\begin{equation}\label{rkw}
\rkw = (\Nsp)^2 \  \in \ \Omega^2(Y)
\end{equation}
determined by
\begin{equation}\label{localrkw}
\rkw_{|U_\s}  \ = \ d\,\omega^{\,\mbox{{\tiny $\Sss(\Ps_\s)$}}}\  \in \ \Omega^2(U_\s).
\end{equation}

\vskip 2mm
{\bf Remark.} \ {\it No use is made of the interpretation of $\Ww$ as a `Frech\'{e}t bundle'. The 2-form $\rkw$ is constructed concretely as the vertical trace of a vertical $\pdo$-valued form on $M$ (cf. Appendix).}

\subsubsection{Why $\rkw$ is a relative eta form}

The APS $\eta$-invariant of a single invertible Dirac-type operator $\pd$ over a closed manifold $N$ is
$$\eta(\pd)  \ = \  \frac{1}{2\sqrt{\pi}}\int^{\o}_0
    t^{-1/2}\Tr(\pd e^{-t\pd^2}) \ dt \ = \ \Tr(\pd|\pd|^{-s-1})|_{s=0}^{{\rm mer}},$$
the superscript indicating the meromorphically continued trace
evaluated at $s=0$. Equivalently,
\begin{equation}\label{e:etareg}
 \eta(\pd) = \Tr((\Pi^{\pd}_{>} - \Pi^{\pd}_{<})|\pd|^{-s})|_{s=0}^{{\rm
 mer}}
\end{equation}
is the zeta function quasi-trace of the involution $\Pi^{\pd}_{>}
- \Pi^{\pd}_{<}$ defined by the order zero $\pdo$ projections
$\Pi^{\pd}_{>} = \frac{1}{2}(I + \pd|\pd|\ii)$ and $\Pi^{\pd}_{<}
= \frac{1}{2}(I - \pd|\pd|\ii) = (\Pi^{\pd}_{>})^\perp$ onto the
positive and negative spectral subspaces of $\pd$.

Consider $\pdo$ projections $P, P^{'}$ with $P - \Pi_{>}$ a and
$P^{'} - \Pi_{>}$ smoothing operators. Since $P - P^{'}$ is
smoothing the relative variant of \eqref{e:etareg} exists without
regularization
\begin{equation}\label{e:releta}
    \eta(P,P^{'}) = \Tr\left((P - P^{\perp}) - (P^{'} -
    (P^{'})^{\perp})\right).
\end{equation}
One then has $\eta(\Pi_{>}^{\pd},\Pi_{>}^{\pd^{'}}) = \eta(\pd) -
\eta(\pd^{'})$ for $\pd -\pd^{'}$ a finite-rank $\pdo$, and the
relative index formula
\begin{equation}\label{e:releta4}
    \frac{\eta(P,P^{'})}{2} = \ind(\pd_P) - \ind(\pd_{P^{'}}),
\end{equation}
which is the pointwise content of \eqref{indexbundle P}. This is
the form degree zero in the boundary Chern character form
$\eta(\Ps^{'},\Ps)$ whose component in $\Omega^{2k}(Y)$  is up to
a constant the vertical trace
\begin{equation*}
\eta(\Ps,\Ps^{'})_{\,[2k]} = \Tr\left((\na^{\,\mbox{{\tiny
${\mathcal W}$}}})^{2k} - (\na^{\,\mbox{{\tiny ${\mathcal
W}$}}^{'}})^{2k}\right).
\end{equation*}
In particular,  $\rkw = \eta(P(\Ds),\Ps)_{\,[2]}$.

\subsection{The zeta function connection on $\Det \Ds_\Ps$ }

The $\z$-connection on $\Det \Ds_\Ps$ is defined locally on $U_\s$ by
\begin{equation}\label{z nabla 1}
\na^{\,\mbox{{\tiny $\zeta, \Pss$}}} \det \Ds_{\Pss} =  \omega^{\,\mbox{{\tiny $\zeta, \Pss$}}}\,\det \Ds_{\Pss}
\end{equation}
with
\begin{equation}\label{zeta U sigma}
\omega^{\,\mbox{{\tiny $\zeta, \Pss$}}} \ = \ -\left.\Tr(\Delta_{P(\Ds)}\si \Ds_{\Pss}\Np \Ds_{\Pss}\ii)\right|^{{\rm mer}}_{s=0} \ \ \in \  \Omega^1(U_\s),
\end{equation}

where $\Tr: \G(Y,\,\Psi^{-\oo}_\flat(E))\  = \ \Psi^{-\oo}_{{\rm vert}, \,\flat}(M,E) \too \Ci(Y)$ is the vertical trace (integral over the fibres, see Appendix).

Here, the notation $\left.\Tr(Q(s)\right|^{{\rm mer}}_{s=0}$ for a family of operators $Q(s)$ depending holomorphically on $s$ and of trace-class for $\re(s)>>0$, means the constant term around $s=0$ (the `finite part') in the Laurent expansion of the meromorphic extension $\left.\Tr(Q(s))\right|^{{\rm mer}}$ of the trace of $Q(s)$ from $\re(s)>>0$ to all of $\C$, assuming this is defined.

The definition of $\omega^{\,\mbox{{\tiny $\zeta, \Pss$}}}$ has particular features which make it work (and be the essentially canonical choice). These are as follows.

The operator $\Ds_{\Pss}$ on the right-side of \eqref{zeta U sigma} means that
\begin{equation}\label{Np in dom Dp}
(\Np_\xi \Ds_{\Pss}\ii) s \in {\rm dom}(\Ds_{\Pss}), \ \ \ \ \ s\in\G(M,E^-).
\end{equation}
Ensuring that \eqref{Np in dom Dp} holds is the job of the connection $\Np$, which is constructed in \secref{Np 1} (this issue is not present in the case of boundaryless manifolds). That is,
\begin{equation}\label{zeta U sigma 2}
\omega^{\,\mbox{{\tiny $\zeta, \Pss$}}}_\s  \ = \ -\,\left.\Tr(\Delta_{P(\Ds)}\si \Ds\Np \Ds_{\Pss}\ii)\right|^{{\rm mer}}_{s=0} \ \ \in \Omega^1(U_\s)
\end{equation}
while the additional subscript in \eqref{zeta U sigma} indicates \eqref{Np in dom Dp}.
$\Np$ has also to be such that the local 1-forms \eqref{zeta U sigma} patch together to define a global $\z$-connection on $\Det \Ds_\Ps$.

\vskip 2mm

The regularized trace of $\Ds_{\Pss}\Np \Ds_{\Pss}\ii$ in \eqref{zeta U sigma} is defined for any vertical APS $\pdo$ projection $\Ps$ using the complex power  $\Delta_{P(\Ds)}\si$ of the Calder\'{o}n Laplacian  $\Delta_{P(\Ds)} = \Ds\Ds^*$ (cf.\eqref{Delta P}).

This differs from the case of boundaryless manifolds which, recall, works as follows. Suppose $\Dd$ is a smooth family of Dirac-type operators associated to a fibration $\pi: N \too Y$ of  compact boundaryless manifolds. Then the determinant line bundle $\Det \Dd$ may be constructed with respect to local charts
$U_{\textsf{s}} = \{y\in Y \ | \ \Dd + \textsf{s} \ {\rm invertible}\}$  with $\textsf{s}\in\Psi^{-\oo}_{{\rm vert}}(N,E^+,E^-)$ a vertical smoothing operator. Over $U_{\textsf{s}}$ one has  the trivialization $y \mtoo \det (\Dd_y + \textsf{s}_y) \in \Det(\Dd_y + \textsf{s}_y)$ and the $\z$-connection 1-form is
$-\left.\Tr\left(\Delta_\textsf{s}\si (\Dd+\textsf{s})\na (\Dd+\textsf{s})\ii\right)\right|^{{\rm mer}}_{s=0}$, where  $\Delta_\textsf{s}$ is the Laplacian of $\Dd+\textsf{s}$.
What makes the patching work in this case is
\begin{equation}\label{patching N}
\left.\Tr(\,(\Delta_\textsf{s}\si - \Delta_{\textsf{s}'}\si)(\Ds+\textsf{s})\na (\Dd+\textsf{s})\ii)\right|^{{\rm mer}}_{s=0}  \ = \ 0 \ \ \ {\rm for } \ \ \textsf{s}, \textsf{s}'\in\Psi^{-\oo}_{{\rm vert}}(N,E^+,E^-).
\end{equation}
This is easily seen; for example, from the precise formulae of \cite{PaSc}. This might suggest that the local $\z$-connection form on $\Det \Ds_\Ps$  be defined as $\left.\Tr\left(\Delta_{\Pss}\si \Ds\Np \Ds_{\Pss}\ii\right)\right|^{{\rm mer}}_{s=0}$. But these forms  do not patch together, because the analogue of the left-side of \eqref{patching N} does not vanish. This one knows from the pole structure of the meromorphic continuation of the trace to all of $\C$, from \cite{GS1, GS2, Gr1, Gr2, Gr3} the constant term in the Laurent expansion at zero depends on $\Pss, \Ps_{\sigma'}$.

In contrast, the connection forms $\omega^{\,\mbox{{\tiny $\zeta, \Pss$}}}$ do patch together  (\thmref{theorem}).

This carries a certain naturality, the family of vertical APS boundary problems $\D_{P(\Ds)}$ is distinguished by the fact that it is invertible (at all points $y\in Y$), and thus so is $\Delta_{P(\Ds)}$, providing a global regularizing operator not available in the case of general family $\Dd$ over boundaryless manifolds. In general, changing the regularizing family $Q(s)$ of elliptic $\pdos$ used to define the connection form
$\left.\Tr\left(Q(s) (\Dd+\textsf{s})\na (\Dd+\textsf{s})\ii\right)\right|^{{\rm mer}}_{s=0}
$ results in additional residue trace terms.

\subsubsection{A connection on $\Hom(\Hh(E^-),\Hh_\Ps(E^+))$}\label{Np 1}

To define a connection on $\Det \Ds_\Ps$ requires a connection on the bundle $\Hom(\,\Hh(E^-),\Hh_\Ps(E^+))$ whose  sections are the subspace of vertical $\pdbos$ with range in $\Ker(\Ps\circ\gam)={\rm dom}(\Ds_\Ps)$
$$\G\left(M,\Hom(\,\Hh(E^-),\Hh_\Ps(E^+))\right) \ :\,=\  \{\As\in \Psi_{{\rm vert},\beta}(M,E^-,E^+) \ | \ \Ps\gam \As s =0, \ s\in \G(M,E^-)\}.$$
All that that requires is a natural connection $\Np$ on $\Hh_\Ps(E^+)$, meaning a connection
$\tNp$ on $\G(Y,\Hh(E^+))=\G(M,E^+)$ which preserves  $\Ker(\Ps\circ\gam)={\rm dom}(\Ds_\Ps).$ That is, such that
\begin{equation}\label{Np on s in dom Dp}
\Ps\gam\Np_\xi s := \Ps\gam\tNp_{\xi_H} s =0 \ \ \ {\rm for}  \ s\in\G(M,E^+) \ \ {\rm with} \ \ \Ps\gam s=0.
\end{equation}
For, then, there is the induced connection (also denoted $\Np$) on $\Hom(\,\Hh(E^-),\Hh_\Ps(E^+))$
$$\Np_\xi  \ : \, \G(M,\Hom(\,\Hh(E^-),\Hh_\Ps(E^+))) \too \G(M,\Hom(\,\Hh(E^-),\Hh_\Ps(E^+))) ,$$
\begin{equation}\label{nabla hom P a}
(\Np_\xi \As) s  :\,= \Np_\xi(\As s) - \As(\Nm_\xi s)  :\,= \tNp_{\xi_H}(\As s) - \As(\wt{\nabla}_{\xi_H} s)
\end{equation}
where $\wt{\nabla}$ is the connection \eqref{Nb} and $\xi\in\G(Y,TY)$. We then evidently have
\begin{equation}\label{nabla hom P b}
\Ps\gam (\Np_\xi \As) s  =  0  \ \ \ \ {\rm for}  \ s\in\G(M,E^-).
\end{equation}
This is how $\Np \Ds_{\Pss}\ii$
in \eqref{zeta U sigma} is defined, and why $\Ds\Np \Ds_{\Pss}\ii$=$\Ds_\Ps\Np \Ds_{\Pss}\ii$.

The task, then, is to define the connection $\Np$ in \eqref{Np on s in dom Dp}.
The connection \eqref{nabla2} $\Nm$ on $\Hh(E^+)$ does not restrict to a connection on  $\Hh_\Ps(E^+)$ (except when $P$ is constant in $y\in Y$ as in the example of \secref{Gr and Diff S1 and CFT}\,), i.e. \eqref{Np on s in dom Dp} does not hold for $\na$. We define $\Np$ by adding a correction term to $\Nm$ in an essentially canonical way, as follows.

First, for an APS-type vertical boundary $\pdo$ $\Ps\in \G(Y,\Psi^0_{{\rm vert},\flat}(\Eb)) :\,= \Psi^0_{{\rm vert}}(\pd M, \Eb)$ we have its covariant derivative
$$ \Nb_\xi \Ps\ \in\  \in \G(Y,\Psi^0_{{\rm vert},\flat}(\Eb)),$$
where $\Nb$ is the connection \eqref{nabla hom}. Let $\phi: [0,\oo) \to \R$ be a smooth function with $\phi(u)=1$ for $0\leq u < 1/4$ and $\phi(u)=0$ for $u>3/4$. Define
$$m_\phi : M \to \R$$
with support in the collar neighbourhood $\Uu$ of $\pd M$ by
$$m_\phi(x) \ = \
\left\{
  \begin{array}{ll}
    0, &    \ \ x\in M\bsh\Uu,\\
    \phi(u), & \ \ x =(u,z) \in \Uu = [0,1)\times \pd M.
  \end{array}
\right.
$$
Then we define
\begin{equation}\label{Np}
\Np\ :=\  \Nm \ +\  m_\phi\, \Ps(\Nb \Ps)\gam.
\end{equation}
Thus for $\xi\in\G(Y,TY)$ and $s\in\G(M,E^+)$
\begin{equation}\label{Np M}
\Np_\xi s\ :\,=\ \tNp_{\xi_H} s\ :\,=\  \tNm_{\xi_H}s \ +\  m_\phi\, \Ps(\tNb_{\xi_H}P)\gam s
\end{equation}
and the second (endomorphism) term acts by
\begin{equation}\label{second term}
(m_\phi\, \Ps(\tNb_{\xi_H}P)\gam s)(x) \ = \
\left\{
  \begin{array}{ll}
    0, &    \ \ x\in M\bsh\Uu,\\
    \phi(u)\, \Ps(\tNb_{\xi_H}P)(s(0,z)) & \ \ x =(u,z) \in \Uu.
  \end{array}
\right.
\end{equation}
Because of the restriction map $\gam$ the Leibnitz property
does not hold for $\Np$ on $\G(M,E^+)$ as a  $\Ci(M)$ module. It does hold, however, for $\G(M,E^+)$ as a $\Ci(Y)$ module \eqref{module}, which is exactly what we need; that is, for the $\Ci(Y)$ multiplication \eqref{module}
\begin{equation}\label{Np leibnitz}
\Np f\cdot s\ =\  df \cdot s + f\cdot \Np s \ \ \ \ {\rm for} \ \ f\in\Ci(Y), \ \ s\in \G(Y,\Hh_\Ps(E^+)).
\end{equation}

\begin{prop}\label{connection on Hp}
$\Np$ defines a connection on $\Hh_\Ps(E^+)$. That is, \eqref{Np on s in dom Dp} holds so that
\begin{equation}\label{Np connection}
\Np_\xi : \G(Y,\Hh_\Ps(E^+)) \too \G(Y,\Hh_\Ps(E^+))
\end{equation}
and satisfies the Leibnitz property \eqref{Np leibnitz}.
\end{prop}

{\bf Proof} \ The Leibnitz property of the first term \eqref{nabla2} of $\Np$ is standard
$$
\Nm_\xi (f\cdot s)  = \wt{\nabla}_{\xi_H} \left((f\circ\pi)s\right)
 = \xi_H\left(f\circ\pi\right).\,s + (f\circ\pi)\wt{\nabla}_{\xi_H}s
 = df(\xi)\,\cdot\,s + f\cdot\Nm_{\xi}s $$
using the Leibnitz property of $\wt{\nabla}$ for the second equality and the chain rule for the third.
Thus \eqref{Np leibnitz} is equivalent to the linearity for $f\in\Ci(Y)$ and $s\in \G(M,E^+)$
\begin{equation}\label{term 2 linear}
m_\phi\, \Ps(\tNb_{\xi_H}P)\gam (f\cdot s)  = f\cdot m_\phi\, \Ps(\tNb_{\xi_H}P)\gam s,
\end{equation}
and this holds because $f$ acts as a constant on each fibre $X_y$ of $M$, by definition \eqref{module}. Precisely, we may assume $ x =(u,z) \in \Uu$, the expressions being zero otherwise, and then from \eqref{second term}
\begin{eqnarray*}
m_\phi\, \Ps(\tNb_{\xi_H}P)\gam (f\cdot s)(u,z) & = & \phi(u)\,\Ps(\tNb_{\xi_H}P)\left(f(\pi(0,z))s(0,z)\right) \\
&= & f(\pi(0,z))\phi(u)\,\Ps(\tNb_{\xi_H}P)(s(0,z))\\ & = &
f(\pi(u,z))\phi(u)\,\Ps(\tNb_{\xi_H}P)(s(0,z))\\
&=& (f\cdot m_\phi\, \Ps(\tNb_{\xi_H}P)\gam s)(u,z)
\end{eqnarray*}
which is \eqref{term 2 linear}. To see \eqref{Np on s in dom Dp}, we have applying $\Ps\circ \gam$ to \eqref{Np M}
\begin{eqnarray}
\Ps\gam \Np_\xi s & = & \Ps\gam\tNm_{\xi_H}s + \Ps(\tNb_{\xi_H}P)(\gam s) \nonumber\\
& = & \Ps\,\tNb_{\xi_H}(\gam s) + \Ps(\tNb_{\xi_H}P)(\gam s),\label{PgamNabla s}
\end{eqnarray}
using \eqref{exact} for the second equality. From \eqref{nabla hom}
\begin{equation*}
\Ps(\nablab_\xi \Ps) (h)  = \Ps\,\tNb_{\xi_H}(\Ps h) - \Ps\tNb_{\xi_H} h, \hskip 10mm h\in \G(\pd M,\Eb).
\end{equation*}
So with $h=\gam s$ and the assumption of \eqref{Np on s in dom Dp}
\begin{equation*}
\Ps(\nablab_\xi \Ps) (\gam s)  = - \,\Ps\,\tNb_{\xi_H}(\gam s), \hskip 10mm h\in \G(\pd M,\Eb),
\end{equation*}
and hence \eqref{PgamNabla s} vanishes.
\begin{flushright}$\square$
\end{flushright}

\subsubsection{Curvature of $\Nzp$ }

The curvature of the connection  $\Nzp$ on the complex line bundle $\Det \Ds_\Ps \too Y$ is the globally defined two form
\begin{equation}\label{fzp}
F_\z(\Ds_\Ps) \ = \ (\Nzp)^2 \  \in \ \Omega^2(Y)
\end{equation}
determined locally by
\begin{equation}\label{FzP}
F_\z(\Ds_\Ps)_{|U_\s}  \ = \ d\,\omega^{\,\mbox{{\tiny $\zeta,\,\Ps_\s$}}}\  \in \ \Omega^2(U_\s).
\end{equation}

\begin{thm}\label{theorem}
The locally defined $\z$ 1-forms \eqref{zeta U sigma} define a connection on the determinant line bundle $\Det \Ds_\Ps$ with curvature
\begin{equation}\label{det P curv}
F_\z(\Ds_\Ps)\  = \ F_\z(\Ds_{\Ps(\Ds)})\  + \ \rkw.
\end{equation}

$F_\z(\Ds_{\Ps(\Ds)})$ is canonically exact, precisely
\begin{equation}\label{beta zeta}
\beta_\z \ :=\ \Tr(\Delta_{P(\Ds)}^z \Ds\,\Nk \Ds_{P(\Ds)}\ii)|^{{\rm mer}}_{z=0} \ \ \in \Omega^1(Y)
\end{equation}
is a globally defined 1-form and
\begin{equation}\label{det P(D) curv}
F_\z(\Ds_{\Ps(\Ds)})\  = \ d\,\beta_\z.
\end{equation}
\end{thm}

\subsection{Proof of \thmref{theorem}}

For the patching of the connection forms, the issue is that there are two candidates for the local connection over $U_\s \cap U_{\sigma'}$ defined by \eqref{z nabla 1}. Let $l$ be a smooth section of $\Det \Ds_\Ps$ over $U_\s \cap U_{\sigma'}$. Then $$l = f_\s.\det \Ds_{\Pss}  = f_{\s'}.\det \Ds_{\Ps_{\s^{'}}}$$
for smooth functions $f_\s, f_{\s'}: U_\s \cap U_{\sigma'}\to \C$. The covariant derivative of $l$
is therefore
$$\na^{\,\mbox{{\tiny $\zeta, \Pss$}}}  (f_\s.\det \Ds_{\Pss})  = df_\s.\det \Ds_{\Pss} + f_\s. \omega^{\,\mbox{{\tiny $\zeta, \Pss$}}}\,\det \Ds_{\Pss} $$
and also
$$\na^{\,\mbox{{\tiny $\zeta, \Ps_{\s'}$}}} (f_{\s'}.\det \Ds_{\,\Ps_{\s^{'}}}) =  df_{\s'}.\det \Ds_{\,\Ps_{\s^{'}}} +  f_{\s'}.\omega^{\,\mbox{{\tiny $\zeta, \Ps_{\mbox{{\tiny $\sigma'$}}}$}}}\,\det \Ds_{\,\Ps_{\s^{'}}}$$
and these must coincide. From \eqref{transition Det Dp}
$$f_{\s'} = {\rm det}_F \left(\Sss(\Pss) \circ \Sss(\Ps_{\mbox{{\tiny $\sigma'$}}})\ii \right)\,f_\s
\hskip 7mm {\rm on} \ \ \ U_\s \cap U_{\sigma'}.$$
Hence, using \secref{connection on Det sp}, the patching condition for the locally defined connection forms is
\begin{equation}\label{connection patching}
\omega^{\,\mbox{{\tiny $\zeta, \Pss$}}} \ -\ \omega^{\,\mbox{{\tiny $\zeta, \Ps_{\mbox{{\tiny $\sigma'$}}}$}}}
\ = \  \omega^{\,\mbox{{\tiny $\Sss(\Ps_\s)$}}}\ - \ \omega^{\,\mbox{{\tiny $\Sss(\Ps_{\mbox{{\tiny $\sigma'$}}})$}}}\hskip 7mm {\rm on} \ \ \ U_\s \cap U_{\sigma'}.
\end{equation}

We will prove a slightly more general statement, which also
captures \eqref{det P curv}. Let $\Ps, \Ps^{'} \in \Psi^0_{{\rm
vert}}(\pd M, \Eb)$ be any two vertical $\pdo$ APS projections and
let $U$ be the open subset of $Y$ where both $(\Ds_{\Ps})_y$ and
$(\Ds_{\Ps^{'}})_y$ are invertible. Then
\begin{equation}\label{rel connection forms}
\omega^{\,\mbox{{\tiny $\zeta, \Ps$}}} \ -\ \omega^{\,\mbox{{\tiny
$\zeta, \Ps^{'}$}}} \ = \ \omega^{\,\mbox{{\tiny $\Sss(\Ps)$}}}\ -
\ \omega^{\,\mbox{{\tiny $\Sss(\Ps^{'})$}}}\hskip 7mm {\rm on} \ \
\ U.
\end{equation}
Here,
\begin{equation*}
\omega^{\,\mbox{{\tiny $\zeta, \Ps$}}} \ = \
-\left.\Tr(\Delta_{P(\Ds)}\si \Ds_{\Ps}\nabla^{\,\mbox{{\tiny
$\textsf{P}$}}} \Ds_{\Ps}\ii)\right|^{{\rm mer}}_{s=0} \hskip 5mm
{\rm and} \hskip 5mm \omega^{\,\mbox{{\tiny $\Sss(\Ps)$}}} =
\Tr(\Sss(\Ps)\ii\nabla_\xi^{\,\mbox{{\tiny ${\mathcal K},{\mathcal
W}$}}}\,\Sss(\Ps)).
\end{equation*}
From \eqref{nabla hom k w}, $\omega^{\,\mbox{{\tiny
$\Sss(P(\Ds))$}}}=0$ and hence, using \eqref{localrkw} and
\eqref{FzP}, \eqref{rel connection forms} also proves \eqref{det P
curv} globally in $\Omega^2(Y)$; note that the right-side of
\eqref{localrkw} and \eqref{FzP} are independent of the choice of
$\s$, i.e. $d\omega^{\,\mbox{{\tiny $\zeta, \Pss$}}}=
d\omega^{\,\mbox{{\tiny $\zeta, \Ps_{\mbox{{\tiny
$\sigma'$}}}$}}}= F_\z(\Ds_\Ps)_{|U_\s \cap U_{\sigma'}}$, and
likewise for $\omega^{\,\mbox{{\tiny $\Sss(\Ps_\s)$}}}$. Clearly,
establishing \eqref{rel connection forms} de facto proves the
identity for the perturbations of $P$ and $P^{'}$ on each chart
$U_\s$, and hence shows \eqref{det P curv} globally.

To see \eqref{rel connection forms}, since the vertical trace
defining the zeta form is taken on $\G(M,E^-)$ (or, rather,
$L^2(M,E^-)$)  we have
\begin{equation}\label{rel zeta forms}
-(\omega^{\,\mbox{{\tiny $\zeta, \Ps$}}}  - \omega^{\,\mbox{{\tiny
$\zeta, \Ps^{'}$}}}) =  \left.\Tr\left(\Delta_{P(\Ds)}\si
\left(\Ds_{\Ps}\nabla^{\,\mbox{{\tiny $\textsf{P}$}}} \Ds_{\Ps}\ii
- \Ds_{\Ps^{'}}\nabla^{\,\mbox{{\tiny $\textsf{P}^{'}$}}}
\Ds_{\Ps^{'}}\ii\right)\right)\right|^{{\rm mer}}_{s=0}.
\end{equation}
(Note that $$\Tr(\Delta_{P(\Ds)}\si
\Ds_{\Ps}\nabla^{\,\mbox{{\tiny $\textsf{P}$}}} \Ds_{\Ps}\ii) -
\Tr(\Delta_{P(\Ds)}\si \Ds_{\Ps^{'}}\nabla^{\,\mbox{{\tiny
$\textsf{P}^{'}$}}} \Ds_{\Ps^{'}}\ii)= \Tr\left(\Delta_{P(\Ds)}\si
\left(\Ds_{\Ps}\nabla^{\,\mbox{{\tiny $\textsf{P}$}}} \Ds_{\Ps}\ii
- \Ds_{\Ps^{'}}\nabla^{\,\mbox{{\tiny $\textsf{P}^{'}$}}}
\Ds_{\Ps^{'}}\ii\right)\right)$$ for large $\Re(s)$, and by the
uniqueness of continuation this extends to all of $\C$.) From
\eqref{P1-P2 a} and \eqref{nabla hom}
\begin{equation}\label{rel Nbp}
\nabla^{\,\mbox{{\tiny $\textsf{P}$}}}  - \nabla^{\,\mbox{{\tiny
$\textsf{P}^{'}$}}}    = m_\phi\, \left(\Ps(\Nb \Ps)- \Ps^{'}(\Nb
\Ps^{'})\right)\gam \ \ \in \ \G(Y, \Psi^{-\oo}_{{\rm
vert,\flat}}(E)),
\end{equation}
and hence  using \propref{relinv smooth}
\begin{equation}\label{rel np Dp }
\nabla^{\,\mbox{{\tiny $\textsf{P}$}}} \Ds_{\Ps}\ii -
\nabla^{\,\mbox{{\tiny $\textsf{P}^{'}$}}} \Ds_{\Ps^{'}}\ii =
\left(\nabla^{\,\mbox{{\tiny $\textsf{P}$}}} -
\nabla^{\,\mbox{{\tiny $\textsf{P}^{'}$}}}\right)\Ds_{\Ps^{'}}\ii
+
 \nabla^{\,\mbox{{\tiny $\textsf{P}$}}}\left(\Ds_{\Ps}\ii -  \Ds_{\Ps^{'}}\ii\right)\ \ \in \ \G(Y, \Psi^{-\oo}_{{\rm vert,\flat}}(E)).
\end{equation}
Hence
$$\Ds_{\Ps}\nabla^{\,\mbox{{\tiny $\textsf{P}$}}} \Ds_{\Ps}\ii - \Ds_{\Ps^{'}}\nabla^{\,\mbox{{\tiny $\textsf{P}^{'}$}}} \Ds_{\Ps^{'}}\ii= \Ds\left(\nabla^{\,\mbox{{\tiny $\textsf{P}$}}} \Ds_{\Ps}\ii - \nabla^{\,\mbox{{\tiny $\textsf{P}^{'}$}}} \Ds_{\Ps^{'}}\ii\right) \ \in \ \G(Y, \Psi^{-\oo}_{{\rm vert,\flat}}(E))$$
is also a smooth family of smoothing operators (with $\Ci$ kernel).  It follows that we may swap the order of the operators inside the trace on the right-side of \eqref{rel zeta forms} to obtain
\begin{eqnarray*}
-(\omega^{\,\mbox{{\tiny $\zeta, \Ps$}}}  - \omega^{\,\mbox{{\tiny
$\zeta, \Ps^{'}$}}})
& = &   \left.\Tr\left( \left(\Ds_{\Ps}\nabla^{\,\mbox{{\tiny $\textsf{P}$}}} \Ds_{\Ps}\ii - \Ds_{\Ps^{'}}\nabla^{\,\mbox{{\tiny $\textsf{P}^{'}$}}} \Ds_{\Ps^{'}}\ii\right)\Delta_{P(\Ds)}\si\right)\right|^{{\rm mer}}_{s=0} \\
& = &
\left.\Tr\left(\Ds\left(\nabla^{\,\mbox{{\tiny $\textsf{P}$}}} \Ds_{\Ps}\ii - \nabla^{\,\mbox{{\tiny $\textsf{P}^{'}$}}} \Ds_{\Ps^{'}}\ii\right)\Delta \Delta_{P(\Ds)}^{-s-1}\right)\right|^{{\rm mer}}_{s=0} \\
& = &
\Tr\left(\Ds\left(\nabla^{\,\mbox{{\tiny $\textsf{P}$}}} \Ds_{\Ps}\ii - \nabla^{\,\mbox{{\tiny $\textsf{P}^{'}$}}} \Ds_{\Ps^{'}}\ii\right)\Delta \Delta_{P(\Ds)}\ii\right) \\
& = &
\Tr\left(\Ds\left(\nabla^{\,\mbox{{\tiny $\textsf{P}$}}} \Ds_{\Ps}\ii - \nabla^{\,\mbox{{\tiny $\textsf{P}^{'}$}}} \Ds_{\Ps^{'}}\ii\right)\right) \\
\end{eqnarray*}
 using that $\Delta_{P(\Ds)}^{-s-1}$ is vertically norm continuous for $\Re(s)>-1$ and, in particular, at $s=0$, and hence that we may take $s$ down to zero without continuation of the vertical trace.

Using \eqref{rel Nbp}, \eqref{rel np Dp } and \eqref{relinv} we therefore have
\begin{eqnarray*}
\omega^{\,\mbox{{\tiny $\zeta, \Ps$}}}  - \omega^{\,\mbox{{\tiny
$\zeta, \Ps^{'}$}}}
= &  &  \Tr\left(\Ds\Nm\left(\Ks(\Ps)\Ps\gam \Ds_{\Ps^{'}} \ii\right)\right) \hskip 65mm {\rm (I)}\\[2mm]
&  & \hskip 10mm -\ \Tr\left(\Ds \,m_\phi\, \left(\Ps(\Nb \Ps)-
\Ps^{'}(\Nb \Ps^{'})\right)\gam \Ds_{\Ps^{'}}\ii \right)
\hskip 24mm {\rm (II)}\\[2mm]
&  & \hskip 20mm  +\ \Tr\left(\Ds \,m_\phi\,\Ps(\Nb
\Ps)P(\Ds)\Sss(\Ps)\ii\Ps\gam \Ds_{\Ps^{'}}\ii \right)\hskip 15mm
{\rm (III)}
\end{eqnarray*}
using the fact that each term is a vertical smoothing operator, as in the proof of \propref{relinv smooth} for terms (I) and (III). We will deal with these terms in reverse order.

\underline{Term (III)}:

Again, in view of \eqref{rel Nbp} we may permute the order of operators in the trace to obtain
\begin{eqnarray*}
{\rm Term \ (III)} & = &
\Tr\left(\Ps\gam \Ds_{\Ps^{'}}\ii \Ds \,m_\phi\,\Ps(\Nb \Ps)P(\Ds)\Sss(\Ps)\ii\Ps\right) \nonumber\\[2mm]
& \stackrel{\eqref{DPmaster}}{=} &
\Tr\left(\Ps\gam(\Is - \Ks(\Ps^{'})\gam)\,m_\phi\,\Ps(\Nb \Ps)P(\Ds)\Sss(\Ps)\ii\Ps\right) \nonumber\\[2mm]
& \stackrel{\eqref{calderon},\ \eqref{Poisson Dp}}{=} &
\Tr\left(\Ps(\Is - P(\Ds)\Sss(\Ps^{'})\ii\Ps^{'})\,\Ps(\Nb \Ps)P(\Ds)\Sss(\Ps)\ii\Ps\right)\nonumber \\[2mm]
& = & \Tr\left(P(\Ds)\left(\Sss(\Ps)\ii\Ps \ - \
\Sss(\Ps^{'})\ii\Ps^{'}\right) \,\Ps(\Nb \Ps)P(\Ds)\right),
\label{term III}
\end{eqnarray*}
circling the operator $P(\Ds)\Sss(\Ps)\ii\Ps =P(\Ds)\circ
P(\Ds)\Sss(\Ps)\ii\Ps$ around for the final equality.

\underline{Term (II)}:

Since  $\Ps\gam \Ds_{\Ps^{'}}\ii =\Ps \circ\Ps\gam(\Is-\Ps^{'})
\Ds_{\Ps^{'}}\ii$ is a composition of vertically smoothing and
$L^2$-bounded operators we may cycle the operators in the trace to
obtain
\begin{eqnarray*}
 {\rm Term \ (II)} & = & -\ \Tr\left(\gam \Ds_{\Ps^{'}}\ii \Ds \,m_\phi\, \left(\Ps(\Nb \Ps)-  \Ps^{'}(\Nb \Ps^{'})\right)\right) \nonumber\\[2mm]
& \stackrel{\eqref{DPmaster}}{=} & -\ \Tr\left(\gam (\Is - \Ks(\Ps^{'})\Ps^{'}) \,\left(\Ps(\Nb \Ps)-  \Ps^{'}(\Nb \Ps^{'})\right)\right) \nonumber\\[2mm]
 & = & -\ \Tr\left(\Ps(\Nb \Ps)-  \Ps^{'}(\Nb \Ps^{'}) - P(\Ds)\Sss(\Ps^{'})\ii\Ps^{'}\left(\Ps(\Nb \Ps)-  \Ps^{'}(\Nb \Ps^{'})\right) \right). \label{term II}
\end{eqnarray*}

\underline{Term (I)}:

From the functoriality of connections on the hom-bundles \eqref{nabla hom}, \eqref{nabla hom P a},
\begin{eqnarray*}
\Nm\left(\Ks(\Ps)\Ps\gam \Ds_{\Ps^{'}} \ii\right) & :\ =&  \Nm\left(\Ks\circ P(\Ds)\Sss(\Ps)\ii\Ps\gam \Ds_{\Ps^{'}} \ii\right)\\[2mm]
& :=&  \Nm(\Ks)\circ P(\Ds)\Sss(\Ps)\ii\Ps\gam \Ds_{\Ps^{'}} \ii \
+ \ \Ks\, \Nb\left(P(\Ds)\Sss(\Ps)\ii\Ps\gam \Ds_{\Ps^{'}}
\ii\right).
\end{eqnarray*}
Hence from \eqref{KerDs=Kk(Ds)}
$$\Ds\Nm\left(\Ks(\Ps)\Ps\gam \Ds_{\Ps^{'}} \ii\right) = \Ds \Nm(\Ks)\circ P(\Ds)\Sss(\Ps)\ii\Ps\gam \Ds_{\Ps^{'}} \ii$$
and therefore
\begin{eqnarray*}
 {\rm Term \ (I)} & = &  \Tr\left(\Ds \Nm(\Ks)\circ P(\Ds)\Sss(\Ps)\ii\Ps\gam \Ds_{\Ps^{'}} \ii\right) \nonumber\\[2mm]
 & = & \Tr\left(P(\Ds)\Sss(\Ps)\ii\Ps\gam \Ds_{\Ps^{'}} \ii \Ds \Nm(\Ks)\,P(\Ds)\right) \nonumber\\[2mm]
 & = & \Tr\left(P(\Ds)\Sss(\Ps)\ii\Ps\gam (\Is - \Ks(\Ps^{'})\Ps^{'}\gam) \Nm(\Ks)\,P(\Ds)\right) \nonumber\\[2mm]
& \stackrel{\eqref{exact}, \ \eqref{calderon}}{=} & \Tr\left(P(\Ds)\Sss(\Ps)\ii\Ps (\Is - \Ks(\Ps^{'})\Ps^{'}) \Nb(P(\Ds))\,P(\Ds)\right) \nonumber\\[2mm]
& = & \Tr\left(P(\Ds)\Sss(\Ps)\ii\Ps (\Is - \Ks(\Ps^{'})\Ps^{'}) \Nb(P(\Ds))\,P(\Ds)\right) \nonumber\\[2mm]
& = & \Tr\left(\,\left(P(\Ds)\Sss(\Ps)\ii\Ps  -
P(\Ds)\Sss(\Ps^{'})\ii\Ps^{'}\right) \Nb(P(\Ds))\,P(\Ds)\right).
\end{eqnarray*}
Summing the expression for terms (I), (II) and (III),
\begin{eqnarray*}
\omega^{\,\mbox{{\tiny $\zeta, \Ps$}}}  - \omega^{\,\mbox{{\tiny
$\zeta, \Ps^{'}$}}}
&=&  \Tr\left(P(\Ds)\,\Sss(\Ps)\ii\Ps(\Nb \Ps)\,P(\Ds) - P(\Ds)\Sss(\Ps^{'})\ii\Ps^{'}(\Nb \Ps^{'})\,P(\Ds)\right. \nonumber\\[2mm]
 & & \hskip 3mm + \ \left. P(\Ds)\,\Sss(\Ps)\ii\Ps(\Nb P(\Ds))\,P(\Ds) - P(\Ds)\Sss(\Ps^{'})\ii\Ps^{'}(\Nb P(\Ds))\,P(\Ds)\right)
 \nonumber\\[2mm]
 & &  \hskip 10mm +\  \Tr\left(\Ps(\Nb \Ps) - \Ps^{'}(\Nb \Ps^{'})\right).
\end{eqnarray*}
From
$$\Sss(\Ps)\ii\nabla^{\,\mbox{{\tiny ${\mathcal K},{\mathcal W}_\sigma$}}}\,\Sss(\Ps) = P(\Ds)\,\Sss(\Ps)\ii\Ps(\Nb \Ps)\,P(\Ds) + P(\Ds)\Sss(\Ps)\ii\Ps(\Nb P(\Ds))\,P(\Ds)$$
we are therefore left with
$$ \omega^{\,\mbox{{\tiny $\zeta, \Ps$}}} \ - \  \omega^{\,\mbox{{\tiny $\zeta, \Ps^{'}$}}}
\ = \ \omega^{\,\mbox{{\tiny $\Sss(\Ps)$}}} \ - \
\omega^{\,\mbox{{\tiny $\Sss(\Ps^{'})$}}}
  \  \  + \ \ \Tr\left(\Ps(\Nb \Ps) - \Ps^{'}(\Nb \Ps^{'})\right).$$
Since $\Ps \in \Psi^0_{{\rm vert}}(\pd M, \Eb)$ is an indempotent
we have
$$\Nb \Ps=\Nb (\Ps^2) = \Ps \Nb (\Ps) + (\Nb \Ps)\circ \Ps$$
and hence composing with $\Ps$ on the left that
$$ \Ps(\Nb \Ps)\circ \Ps = 0.$$
Hence
$$\Tr\left(\Ps(\Nb \Ps) - \Ps^{'}(\Nb \Ps^{'})\right) = \Tr\left(\Ps
(\Nb \Ps) \circ \Ps^{\perp} - \Ps^{'}(\Nb \Ps^{'})\circ
(\Ps^{'})^{\perp}\right)$$  with  $\Ps^{\perp} = \Is - \Ps$.
Writing  the operator inside the trace as $$\Ps \Nb \Ps \circ
(\Ps^{\perp} - (\Ps^{'})^{\perp}) +  (\Ps \Nb \Ps - \Ps^{'} \Nb
\Ps^{'})\circ (\Ps^{'})^{\perp},$$ the bracketed vertical $\pdos$
are smoothing and we may cycle the operators through the trace
leaving
$$\Tr\left(\Ps(\Nb \Ps) - \Ps^{'}(\Nb \Ps^{'})\right)\ = \
\Tr\left(\Ps^{\perp}\Ps \Nb \Ps \circ \Ps^{\perp} -
(\Ps^{'})^{\perp}\Ps^{'}(\Nb \Ps^{'})\circ
(\Ps^{'})^{\perp}\right)=0.$$
\begin{flushright}
$\Box$
\end{flushright}

\appendix

{\Large {\bf Appendix}}

\section{Vertical Pseudodifferential Operators}\label{pdo closed fibrations}

\subsection{Families of  $\pdos$ on a closed manifold}\label{families pdos}

A smooth family of  $\pdos$ of constant order $\mu$ associated to
a fibration $\pi: N \stackrel{X}{\too} Y$ of compact boundaryless
manifolds, with $\dim(X)=n$,  with vector bundle $E\to N$ means a
classical $\pdo$ $$\As :\G(N,E^+)\too\G(N,E^-)$$ with Schwartz
kernel $k_\As\in\Dd^{'}(N\times_{\pi}N, E\boxtimes E)$ a vertical
distribution, where the fibre product $N\times_{\pi}N$ consists of
pairs $(x,x^{'})\in N\times N$ which lie in the same fibre, i.e.
$\pi(x) = \pi(x^{'})$,  such that in any local  trivialization
$k_\As$ is an oscillatory integral with vertical symbol $\as\in
S^{\nu}_{{\rm vert}}(N/Y)$ of order $\nu$. Here, $\xi$ is {\em
restricted} to the vertical momentum space, along the fibre. We
refer to $\As$ as a vertical $\pdo$ associated to the fibration
and denote this subalgebra of $\pdos$ on $N$ by
$$\G(Y,\Psi^\nu(E^+,E^-)) = \Psi^\nu_{{\rm vert}}(N,E^+,E^-).$$

Thus in any local trivialization $N_{|U_Y} \cong U_Y\times X_{y}$
over an open subset $U_Y\subset Y$ with $y\in U_Y$, and a
trivialization $E \cong U_Y \times V_y \times \R^N$ of $E$ with
$V_y$ an open subset of $X_y$, a vertical amplitude of constant
(in $y$) order $\nu$ is an element
\begin{equation}\label{e:verticalsymbol}
\as =\as(y,x,x^{'},\xi)\in\G\left( \, U_Y \times (V_y\times
V_y)\times\R^n\backslash\{0\},  \ \End(\R^N)\right)
\end{equation}
satisfying the estimate on compact subsets $K\<N$
\begin{equation}\label{e:verticalschwartzclass}
|\pd^{\alp}_x\pd^{\gam}_{x^{'}}\pd^{\del}_y\pd^{\bet}_{\xi}\as| <
C_{\alp,\gam,\del,\bet,K}(1 + |\xi|)^{\nu - |\bet|}.
\end{equation}
We denote this as $\as\in S^{\nu}_{{\rm vert}}(N/Y)$. Here, $\xi$
may be identified with an element of the vertical (or fibre)
cotangent space $T^{*}_x(N/Y)$. The kernel of $\As$ is then {\em
locally} written on $U_Y \times V_y$ as the distribution with
singular support along the diagonal
\begin{equation*}
k_\As(y,x,x^{'}) = \int_{\R^n} e^{i(x-x^{'}).\,\xi} \
\as(y,x,x^{'},\xi) \ \db\xi.
\end{equation*}
If $\As$ has order $\nu < -n$ this integral is convergent and with
respect to a vertical volume form $d_{\mbox{{\tiny N/Y}}}x$ the
trace $\Tr \As$ is the smooth function on $Y$
\begin{equation*}
(\Tr \As)(y)  = \int_{N/Y}\tr(k_\As(y,x,x))  \ d_{\mbox{{\tiny
N/Y}}}x \ \ \in \Ci(Y), \hskip 10mm \nu < -n.
\end{equation*}

If   $\ws\in S^{-\oo}_{{\rm vert}}(N/Y) = \bigcap_\nu
S^{\nu}_{{\rm vert}}(N/Y)$ then the kernel is an element
$$k_\Ws\in \G(N\times_\pi N, E^*\boxtimes E)$$
and defines a vertical smoothing operator (smooth family of smoothing operators)
$$\Ws\in \G(Y,\Psi^{-\oo}(E)).$$
Any vertical $\pdo$ of order $\nu$ may be written in the form $\As
= {\rm OP}(\as) + \Ws$ with $\as = \as(y,x,\xi)$ a vertical symbol
and $\Ws$ a vertical smoothing operator. Assuming this
representation, we will consider here only classical vertical
$\pdo$s, meaning that the symbol has an asymptotic expansion
$\as\sim\sum_{j\geq 0}\as_j$ with $\as_j$ positively homogeneous
in $\xi$ of degree $\nu-j$. The leading symbol $\as_0$ has an
invariant realization as a smooth section
\begin{equation*}
\as_{0}\in \G(T^*(N/Y),\varphi^*(\End(E))) \ ,
\end{equation*}
where $\varphi : T^*(N/Y) \too N$. If  $\as_0$ is an invertible
bundle map then $\As\in \G(Y,\Psi(E))$ is said to be an {\it
elliptic family}. If there  exists $\theta$  such that $\as_0 -
\lam\Ii$ is invertible for each $\lam\in R_{\theta} =
\{re^{i\theta} \ | \ r > 0\}$, where $\Ii$ is the identity bundle
operator, then $\As$ is elliptic with principal angle $\theta$. In
the latter case one has the resolvent family
\begin{equation}\label{res family}
(\As -\lam\Is)\ii \ \in\G(Y,\Psi^{-\nu}(E))
\end{equation}
and the complex powers
\begin{equation}\label{cx powers}
\As_\theta^z := \frac{i}{2\pi}\int_\Cc \lam^z_\theta\,(\As -\lam\Is)\ii \ d\lam \ \in\G(Y,\Psi^{z\nu}(E)),
\end{equation}
where $\Cc$ is a contour running in along $R_{\theta}$ from infinity to a small circle around the origin, clockwise around the circle, then back out to infinity along $R_{\theta}$, as accounted for in detail in \cite{Sc2}. A principal angle, and hence the complex powers, can only exist if the pointwise index is zero.

For example, if $\ps \in S^{m}_{{\rm vert}}(N/Y)$ is a polynomial of order $m\in\Nline$ in $\xi$ and elliptic, then the corresponding vertical $\pdo$ $\Ds \in  \G(Y,\Psi^m(E^+, E^-))$ is a smooth family of elliptic differential operators of order $m$. Specifically, this is the case for a geometric fibration of Riemannian spin manifolds \secref{spinor} with associated smooth family of twisted Dirac operators \cite{B,BGV}. (The space  $\G(Y,\Psi(E^+, E^-))$ of vertical $\pdo$s between different bundles $E^\pm$ is defined by a trivial elaboration of the above.)

\subsection{Families of pseudodifferential boundary operators}\label{families pdbos}

Let $\pi : M \too Y$ be a smooth fibration of compact manifolds
with boundary with vector bundle $E\to M$ and let $$\wt{\pi} :
\wt{M}\too Y$$  be a smooth fibration of compact boundaryless
manifolds with vector bundle $\wt{E}\to M$ such that $$M\< \wt{M},
\ \ \ \wt{E}_{|M} = E,  \ \ \  {\rm and} \ \ \ \wt{\pi}_{|M} =
\pi.$$ We consider the following vertical families of
pseudodifferential boundary operators (vertical $\pdbo$s) as
defined in the single operator case  by Grubb \cite{Gr2},
elaborating the algebra of Boutet de Monvel. First, one has the
truncated, or restricted, $\pdos$
\begin{equation}\label{A+}
\As_+ : \G(M,E) \too \G(M,E), \hskip10mm \As_+ = r^+\As\e^+,
\end{equation}
where
$$ \As \in \Psi_{{\rm vert}}(\wt{M},\wt{E})$$
is a vertical $\pdo$ associated to the fibration of closed manifolds, and
$$ r^+ : \G(\wt{M},\wt{E}) \too \G(M,E),  \ \ \ \  e^+ : \G(M\bsh \pd M ,E) \too \G(\wt{M},\wt{E}) $$
are the `brutal' restriction and extension-by-zero operators. To avoid $e^+$ or $r^+$ introducing any new singularities $\As$ is assumed to be of integer order and to satisfy the {\it transmission condition} at $\pd M$, which in local coordinates $(x_n,w)\in\Uu =(0,1]\times \pd M$ in the collar neighbourhood of the boundary of $M$ is the requirement
$$\pd^\bet_x\pd^\alp_\xi \as_j(0,w,-\xi_n,0) = (-1)^{\nu-j-|\alp|}\pd^\bet_x\pd^\alp_\xi \as_j(0,w,\xi_n,0) \ \ {\rm for} \ |\xi_n| \geq 1,$$
where $\xi = (\xi_n,\xi_w)$ relative to \eqref{TUu}.

More generally, one considers
\begin{equation}\label{A plus sgo}
\As_+  + \Gs : \G(M,E) \too \G(M,E)
\end{equation}
with $\Gs$  a vertical singular Green's operator (vertical sgo), which we return to in a moment.

A {\it vertical trace operator} of order $\mu\in\R$ and class $r\in\Nline$ is an operator from interior to boundary sections of the form
\begin{equation}\label{trace operators}
\Ts : \G(M,E) \too \G(\pd M,\Eb), \hskip10mm \Ts = \sum_{0\leq j<r}\Sss_j\gam_j + \Ts^{'},
\end{equation}
where the $\Sss_j\in\Psi_{{\rm vert}}(\pd M, \Eb)$ are vertical $\pdos$ on the boundary fibration of closed manifolds, as in \secref{pdo closed fibrations}, while $\gam_j s(x_n,w) = \pd^j_{x_n}s(0,w)$ are the restriction maps to the boundary. The additional term is an operator of the form  $\Ts^{'} = \gam \As_+$ for some restricted $\pdo$ \eqref{A+}.

 A {\it vertical Poisson operator} of order $\mu\in\R$ here will be  an operator from boundary to interior sections of the form
\begin{equation}\label{poisson operators}
\Ks : \G(\pd M,\Eb) \too \G(M,E), \hskip10mm \Ks =  r^+ \Bs\gam^*\Cc, \hskip 10mm \gam : \G(\wt{M},\wt{E}) \to \G(\pd M, \Eb),
\end{equation}
where $\Bs \in \Psi^{\mu-1}_{{\rm vert}}(\wt{M},\wt{E})$ is a family of $\pdos$ in the sense of \secref{pdo closed fibrations}, while $\Cc\in\Psi^m_{{\rm vert}}(\pd M, \Eb)$ is a  vertical differential operator on the boundary fibration of order $m$; however, $m=0$ in the following. Note, the restriction map $\gam$ is here coming from $\wt{M}$, rather than $M$ (same notation).

To make the composition rules work one includes vertical sgo operators in \eqref{A plus sgo} of order $\nu$ and class $r\in\Nline$, these have the form
\begin{equation}\label{sgos}
\Gs = \sum_{0\leq j<r}\Ks_j\gam_j + \Gs^{'},
\end{equation}
where $\Ks_j$ is a vertical Poisson operator of order  $\nu-j$, and  $\Gs^{'}$ is defined in local coordinates near $\pd M$ by an oscillatory integral on a sgo symbol $\gs$ satisfying standard estimates in $\xi$ \cite{Gr2}.

As with closed manifolds if $\As + \Gs $  in \eqref{A plus sgo} has order $\nu < -n$ and we assume the order of $\Cc$ in \eqref{poisson operators} is $m=0$ then the ` distribution kernel' is continuous and the trace $\Tr \As$ is the smooth function, or differential form for de Rham valued symbols, on $Y$
\begin{equation}\label{trace sgo}
\Tr (\As+\Gs)(y)  = \int_{M/Y}k_{\,\As+\Gs}(y,x,x)  \, d_{\mbox{{\tiny M/Y}}}x \ \ \in \Aa(Y), \hskip 10mm ({\rm order}(\As + \Gs) << 0).
\end{equation}

For each of the above classes of $\pdbos$  one considers the subclass of operators defined by polyhomogeneous symbols, appropriately formulated \cite{Gr2}. We denote the resulting algebra by
$$\G(Y,\,\Psi_\flat(E))\  = \ \Psi_{{\rm vert}, \,\flat}(M,E).$$
When the kernel is an element
$$k_{\,\As+\Gs}\in \G(M\times_\pi M, E^*\boxtimes E)$$
then the operator defines a vertical smoothing operator (smooth family of smoothing operators)
$$\As+\Gs \in \G(Y,\,\Psi^{-\oo}_\flat(E))\  = \ \Psi^{-\oo}_{{\rm vert}, \,\flat}(M,E).$$

We refer to \cite{Gr2} and references therein for a precise account of the pseudodifferential boundary operator calculus, which extends to the case of vertical operators in a similar way to the case for compact boundaryless manifolds \cite{Sc2}.

\vskip 8mm

{\small \noi \textsc{Department of Mathematics, King's College
London.}}


\begin{thebibliography}{99}

\bibitem{BGV}  Berline, N., E. Getzler, M. Vergne: `Heat
Kernels and Dirac Operators'.  Grundlehren der Mathematischen
Wissenschaften {\bf 298}, Springer-Verlag, Berlin, 1992.

\bibitem{B} Bismut, J-M.: 1986, `The Atiyah-singer index theorem for
families of Dirac operators: Two heat equation proofs',
Invent. Math {\bf 83}, 91--151.

\bibitem{BiCh} Bismut, J-M., Cheeger, J.: 1990,
`Families index for manifolds with boundary, superconnections
and cones' I and II,  J. Funct. Anal. {\bf 89}, 313--363, and
{\bf 90}, 306--354.

\bibitem{BF} Bismut, J-M., Freed, D.S.: 1986,
`The analysis of elliptic families I', Comm. Math. Phys. {\bf
106}, 159--176.

\bibitem{BL} Bruening, J., Lesch, M.:1999, `On the eta-invariant of certain
non-local boundary value problems', Duke Math. J. {\bf 96},
425-468.

\bibitem{BiSo} Birman M.S., Solomyak  M.Z.: 1977,
`Estimates of singular numbers of integral operators', Russ. Math. Surveys. {\bf
32}, 15--89.

\bibitem{BoWo} Boo\ss--Bavnbek, B., Wojciechowski, K.P.:
1993, `Elliptic Boundary Problems for Dirac Operators',
Birkh\"auser, Boston.

\bibitem{Ca} Calder\'{o}n, A.P.: 1976, `Lecture Notes on Pseudo-differential Operators and Elliptic Boundary Value Problems'. I. Consejo Nacional de Investigaciones Cient\'{i}ficas y T\'{e}cnicas, Inst. Argentino de Matem\'{a}tica, buenos Aires.

\bibitem{Gr1} Grubb, G.: 1999, `Trace expansions for
pseudodifferential boundary problems for Dirac-type operators and
more general systems', Ark. Mat. {\bf 37}, 45--86.

\bibitem{Gr2} Grubb, G.: 2001, `A weakly polyhomogeneous calculus
for pseudodifferential boundary problems',  J. Funct. Anal.
{\bf 184}, 19-76.

\bibitem{Gr3} G. Grubb, `A resolvent approach to
traces and zeta Laurent expansions', AMS Contemp. Math. Proc.,
{\bf 366}, 67--93 (2005). Also  arXiv: math.AP/0311081.

\bibitem{GS1} Grubb, G., Seeley, R.: 1995, `Weakly parametric
pseudodifferential operators and Atiyah-Patodi-Singer boundary
problems',  Invent. Math. {\bf 121}, 481-529.

\bibitem{GS2} Grubb, G., Seeley. R.: 1996,
`Zeta and eta functions for Atiyah-Patodi-Singer operators',
J. Geom. Anal. {\bf 6}, 31--77.

\bibitem{MePi} Melrose, R.B., Piazza, P.: 1997, `Families of Dirac operators,
boundaries and the $b$-calculus`, J. Diff. Geom. {\bf 46}, 99--167.

\bibitem{Mo} Moriyoshi, H.: 1994, `The Euler and Godbillon-Vey forms and symplectic structures on
$Diff^\oo_+(S^1)/SO(2)$`, Contemp. Math. {\bf 179}, 193--203.

\bibitem{P} Piazza, P.: 1996, `Determinant bundles, manifolds with boundary and surgery I',
 Comm. Math. Phys. {\bf 178}, 597--626.

\bibitem{PaSc} Paycha, S., Scott, S.: 2007,
'A Laurent expansion for regularized integrals of holomorphic symbols',
Geom. and Funct. Anal. {\bf 17}, 491 - 536.

\bibitem{Qu} Quillen, D.G.: 1985, `Determinants of Cauchy-Riemann operators over a
Riemann surface', Funk. Anal. i ego Prilozhenya {\bf 19}, 37--41.

\bibitem{Qu2} Quillen, D.G.: 1985, `Superconnections and the Chern character',
 Topology {\bf 24}, 89--95.

\bibitem{Sc3} Scott, S.G.: 1995, `Determinants of Dirac boundary value
problems over odd--dimensional manifolds', Comm. Math. Phys. {\bf 173}, 43--76.

\bibitem{Sc1} Scott, S.: 2002, `Zeta determinants on manifolds with boundary',
 Jour. Funct. Anal. {\bf 192}, 112-185.

\bibitem{Sc2} Scott, S.: 2007, `Zeta forms and the local family index theorem',
Trans. Am. Math. Soc. {\bf 359}, 1925 -1957.

\bibitem{Scott book} Scott, S.:  `Traces and Determinants of Pseudodifferential Operators',
OUP, Math. Monographs, to appear 2008.

\bibitem{ScWo} Scott, S., Wojciechowski, K.,: 2000, 'The zeta-determinant and Quillen determinant over odd dimensional manifolds', Geom. and Funct. Anal. {\bf 10}, 1202-1236.

\bibitem{S1} Seeley, R. T.: 1966, `Singular integrals and boundary value problems',
Amer. J. Math. {\bf 88}, 781--809.

\bibitem{S2} Seeley, R. T.: 1969,  `Topics in pseudodifferential operators'.
In: CIME Conference on Pseudo-Differential operators (Stresa 1968), pp. 167--305. Cremonese.

\bibitem{Se} Segal, G.: 2000, `Lectures on QFT', Stanford Lectures.

\end{thebibliography}
\end{document}